\documentclass[11pt,a4paper]{article}
\usepackage{amsmath}
\usepackage{amsfonts}
\usepackage{amssymb}
\usepackage{amsthm}
\usepackage{mathrsfs}
\usepackage{dsfont}
\usepackage{paralist}
\usepackage{natbib}
\usepackage{bm}
\usepackage{bbm}
\usepackage{color}
\usepackage{graphicx}
\usepackage{tabularx}
\usepackage[left=3.5cm,right=3.5cm,top=2cm,bottom=2cm,includeheadfoot]{geometry}
\usepackage{comment}
\usepackage{tikz}
\usetikzlibrary{decorations.pathreplacing}
\usetikzlibrary{decorations.pathmorphing}
\usepackage{resizegather}
\DeclareFontFamily{OT1}{pzc}{}
\DeclareFontShape{OT1}{pzc}{m}{it}{<-> s * [1.10] pzcmi7t}{}
\DeclareMathAlphabet{\mathpzc}{OT1}{pzc}{m}{it}
\usepackage[colorlinks=true,linkcolor=blue,citecolor=blue,pdfborder={0 0 0}]{hyperref}

\newtheoremstyle{normal}
{2ex}               
{3ex}               
{}                  
{}                  
{\bfseries} 
{}                  
{2pt}   
{\thmname{#1}\thmnumber{ #2.} \thmnote{(#3)}}

\newtheoremstyle{italic}
{2ex}
{3ex}
{\itshape}
{}
{\bfseries} 
{}
{2pt}
{\thmname{#1}\thmnumber{ #2.} \thmnote{(#3)}}

\theoremstyle{normal}
\newtheorem{definition}{Definition}[section]
\newtheorem{remark}[definition]{Remark}

\newtheorem{condition}[definition]{Condition}

\theoremstyle{italic}
\newtheorem{theorem}[definition]{Theorem}
\newtheorem{lemma}[definition]{Lemma}

\newtheorem{corollary}[definition]{Corollary}

\renewcommand{\P}{\mathbb{P}}

\newcommand{\E}{\mathbb{E}}

\newcommand{\1}{\mathbbm{1}}
\newcommand{\Nn}{N_n(1)}
\DeclareMathAccent{\verywidehat}{\mathord}{largesymbols}{'144}

\newcommand{\N}{\mathbb{N}}
\newcommand{\R}{\mathbb{R}}

\newcommand{\al}{\alpha}
\newcommand{\be}{\beta}
\newcommand{\la}{\lambda}

\newcommand{\ka}{\kappa}

\newcommand{\si}{\sigma}

\newcommand{\eps}{\varepsilon}

\newcommand{\De}{\Delta}

\newcommand{\gG}{\mathcal G}

\newcommand{\fF}{\mathcal F}

\newcommand{\lL}{\mathcal L}

\newcommand{\tols}{~\stackrel{\lL-(s)}{\longrightarrow}~}

\newcommand{\pn}{\stackrel{\P}{\longrightarrow}}

\newcommand{\Ind}{\mathds{1}}

\newcommand{\mup}{\mu_{p,\beta}}
\newcommand{\sigs}{|\overline{\sigma\lambda}|^p}
\newcommand{\sigl}{|\overline{\sigma}|^p}
\newcommand{\Ci}{{\mathcal{C}_i^n}}

\allowdisplaybreaks[1]

\begin{document}

\title{On the estimation of the jump activity index in the case of random observation times}

\author{Adrian Theopold and Mathias Vetter\thanks{Christian-Albrechts-Universit\"at zu Kiel, Mathematisches Seminar, Heinrich-Hecht-Platz\ 6, 24118 Kiel, Germany.
{E-mail:} theopold@math.uni-kiel.de and vetter@math.uni-kiel.de.} \bigskip \\
{Christian-Albrechts-Universit\"at zu Kiel}
}

\maketitle

\begin{abstract}
We propose a nonparametric estimator of the jump activity index $\be$ of a pure-jump semimartingale $X$ driven by a $\be$-stable process when the underlying observations are coming from a high-frequency setting at irregular times. The proposed estimator is based on an empirical characteristic function using rescaled increments of $X$, with a limit which depends in a complicated way on $\be$ and the distribution of the sampling scheme. Utilising an asymptotic expansion we derive a consistent estimator for $\be$ and prove an associated central limit theorem. 
\end{abstract}

\medskip



\section{Introduction}
\def\theequation{1.\arabic{equation}}
\setcounter{equation}{0}

Recent years have seen a notable development in the statistical analysis of time-continuous stochastic processes beyond the somewhat classical case of an It\^o semimartingale driven by a Brownian motion. Generalisations of that class of processes are in fact manifold, and one can mention for example the analysis of integrals with respect to fractional Brownian motion (e.g.\ \cite{broufuka2018}, \cite{bibinger2020}), the discussion of L\'evy-driven moving averages (e.g.\ \cite{basetal2017}, \cite{basetal2018}), inference on the solution of stochastic PDEs (e.g.\ \cite{bibitrab2020}, \cite{chong2020}, \cite{kainuchi2021}) and the behaviour of integrals with respect to stable pure-jump processes (e.g.\ \cite{heinpodo2021}, \cite{todorov2015}). All of the above results are concerned with high-frequency observations of the respective processes, always in the case of regularly spaced observations in time. 

On the other hand, it is well understood that the underlying assumption of a regular spacing constitutes an ideal setting that simplifies the theoretical statistical analysis but is typically not met in practical applications. For this reason there has always been a lot of interest in understanding the impact of irregular sampling schemes on the proposed statistical methods. For semimartingales driven by Brownian motion one can mention \cite{hayetal2011} and \cite{mykzha2012} among others, with an almost complete treatment in Chapter 14 of \cite{discret}. The case of Brownian semimartingales with jumps is treated for example in \cite{bibvet2015} and \cite{marvet2019}. All of the aforementioned papers deal with exogenous sampling schemes, i.e.\ when the observation times are essentially independent from the underlying processes. There is also limited research on endogenous observations times, mostly modelled via hitting times. See for example \cite{fukros2012} in the continuous case or \cite{vetzwi2017} when additional jumps are present. 

In this paper we are discussing the case of a near stable jump semimartingale observed at irregular times, i.e.\ the underlying process is given by
\begin{align}\label{eq:Xrep1}
X_t=X_0+\int_{0}^{t}\alpha_s ds+\int_{0}^{t}\sigma_{s-}dL_s+Y_t, \quad t>0, 
\end{align}
and the observation times follow a version of the restricted discretisation scheme from \cite{discret}, essentially providing observation times independent of $X$. Loosely speaking and made precise below, $L$ is driven by a $\be$-stable process and $Y$ comprises the residual jumps while $\al$ and $\si$ are appropriately chosen adapted processes. Our goal in this work is to provide a consistent estimator for $\be$ and to establish an associated central limit theorem. Statistical inference on $\be$ has already been conducted in \cite{todorov2015} and \cite{todorov2017} for regular observations while \cite{Jacod2018} provides the theory in a general model where microstructure noise is present and dominates the statistical analysis. 

A first glance, our strategy to estimate $\be$ somewhat resembles the procedure from \cite{todorov2015}, but there are some notable challenges that might occur in other situations as well. First, we are computing an empirical characteristic function $\widetilde{L}^n(p,u)$ which is constructed from local increments of $X$ (and with an auxiliary parameter $p$), but it is important here to rescale any of these increments relative to the length of the underlying time period. Secondly, one can show convergence of this empirical distribution function to a function $L(p,u,\be)$ which not only is a function of $u$, $p$ and the unknown $\be$ but also also depends specifically on the distribution of the discretisation scheme. Unlike in \cite{todorov2015}, where a consistent estimator for $\be$ is obtained via a suitable functional of empirical characteristic functions computed at arbitrary values $u$ and $v$, we have to use sequences $u_n$ and $v_n$ converging to zero plus an asymptotic expansion to obtain a consistent estimator. This procedure also leads to a drop in the rate of convergence in the associated central limit theorem. 

The remainder of this work is as follows: Section \ref{sec:setting} deals with the assumptions on $X$ as well as on the discretisation scheme. In Section \ref{sec:results} we establish our statistical method and we also present the main results on the asymptotic properties both of $\widetilde{L}^n(p,u)$ and of $\hat{\beta}(p,u_n,v_n)$. A thorough simulation is study is provided in Section \ref{sec:simul} where we also discuss issues connected with the estimation of the asymptotic variance in the normal approximation. All proofs are gathered in Section \ref{sec:proof}.

\section{Setting} \label{sec:setting}
\def\theequation{2.\arabic{equation}}
\setcounter{equation}{0}

Throughout this work we adopt the setting from \cite{todorov2015} and assume that we are given a univariate pure-jump semimartingale as defined in \eqref{eq:Xrep1}, i.e.\ that we observe
\begin{align*}
X_t=X_0+\int_{0}^{t}\alpha_s ds+\int_{0}^{t}\sigma_{s-}dL_s+Y_t
\end{align*}
where $L$ and $Y$ are pure-jump It\^o semimartingales and $\al$ and $\si$ are c\`adl\`ag. All processes are defined on some filtered probability space $(\Omega,\mathcal{F},(\mathcal{F}_t)_{t\geq 0},\mathbb{P})$.

Specific assumptions on these processes will be given below, and we start with conditions on the jump processes $L$ and $Y$. Below, $\kappa(x)$ denotes a truncation function, i.e.\ it is the identity in a neighbourhood around zero, odd, bounded and equals zero for large values of $|x|$. We also set $\kappa'(x)=x-\kappa(x)$, and whenever we discuss the characteristic triplet of a L\'evy process it is to be understood as with respect to this choice of the truncation function.

\begin{condition}\label{ass:general}	
We  impose the following conditions on the processes $L$ and $Y$: 
	\begin{enumerate}
		\item[(a)] $L$ is a L\'evy process with characteristic triplet $(0,0,F)$ where the Lebesgue density of the L\'evy measure $F(dx)$ is given by 
		\begin{align*}
		h(x)=\frac{A}{|x|^{1+\beta}}+\tilde{h}(x) 
		\end{align*}
		for some $\be \in (1,2)$ and some $A > 0$. The function $\tilde{h}(x)$ satisfies $$|\tilde{h}(x)|\leq \frac{C}{|x|^{1+\beta'}}$$ for some $\be' < 1$ and all $|x|\leq x_0$, for some $x_0 > 0$. 
		\item[(b)] $Y$ is a finite variation jump process of the form
		\[
		Y_t=\int_0^t \int_{\R} x \mu^Y(ds,dx)
		\]
		where $\mu^Y(ds,dx)$ denotes the jump measure of $Y$ and its compensator is given by $ds \otimes\nu_s^Y(dx)$. The process $$\left(\int_{\R}\left(|x|^{\beta'}\wedge 1\right)\nu_t^Y(dx)\right)_{t\geq 0}$$ is locally bounded for the parameter $\beta'$ from (a).  
	\end{enumerate}	
\end{condition}

Condition \ref{ass:general} should be read in such a way that the pure-jump L\'evy process $L$ is essentially $\be$-stable while all other jumps (both in $L$ and in $Y$) are of much smaller activity and will be dominated by the $\be$-stable part at high frequency. Note that dependence between $L$ and $Y$ is possible, and this will hold for the jump parts of $\alpha$ and $\sigma$ as well. 

\begin{condition}\label{ass:sigma}
	The processes $\alpha$ and $\sigma$ are It\^o semimartingales of the form
	\begin{align*}
	\alpha_t=\alpha_0&+\int_{0}^{t}b_s^\alpha ds+\int_{0}^{t}\eta_s^\alpha dW_s+\int_{0}^{t}\widetilde{\eta}_s^\alpha d\widetilde{W}_s+\int_{0}^{t}\overline{\eta}_s^\alpha d\overline{W}_s\\&+\int_{0}^{t}\int_E\kappa(\delta^\alpha(s,x))\underline{\widetilde{\mu}}(ds,dx)+\int_{0}^{t}\int_E\kappa'(\delta^\alpha(s,x))\underline{\mu}(ds,dx),\\
	\sigma_t=\sigma_0&+\int_{0}^{t}b_s^\sigma ds+\int_{0}^{t}\eta_s^\sigma dW_s+\int_{0}^{t}\widetilde{\eta}_s^\si d\widetilde{W}_s+\int_{0}^{t}\overline{\eta}_s^\si d\overline{W}_s\\&+\int_{0}^{t}\int_E\kappa(\delta^\sigma(s,x))\underline{\widetilde{\mu}}(ds,dx)
	+\int_{0}^{t}\int_E\kappa'(\delta^\sigma(s,x))\underline{\mu}(ds,dx)
	\end{align*}
	where
	\begin{itemize}
		\item[(a)] $|\sigma_t|$ and $|\sigma_{t-}|$ are strictly positive;
		\item[(b)] $W$, $\widetilde{W}$ and $\overline W$ are independent Brownian motions, $\underline{\mu}$ is a Poisson random measure on $\R_+\times E$  with compensator $dt\otimes\lambda(dx)$ for some $\sigma$-finite measure $\lambda$ on $E$ and $\underline{\widetilde{\mu}}$ is the compensated jump measure;
		\item[(c)] $\delta^\alpha(t,x)$ and $\delta^\sigma(t,x)$ are predictable with $|\delta^\alpha(t,x)|+|\delta^\sigma(t,x)|\leq \gamma_k(x)$ for all $t\leq T_k$, where $\gamma_k(x)$ is a deterministic function on $\R$ with $\int_E\left(|\gamma_k(x)|^r\wedge 1\right)\lambda(dx)<\infty$ for some $0\leq r < 2 $ and $T_k$ is a sequence of stopping times increasing to $+\infty$;
		\item[(d)] $b^\alpha, b^\sigma$ are locally bounded while $\eta^\alpha, \eta^\sigma, \widetilde{\eta}^\alpha , \widetilde{\eta}^\si, \overline{\eta}^\alpha $ and $\overline{\eta}^\si$ are c\`adl\`ag. 
	\end{itemize}	
\end{condition}

These assumptions on $\al$ and $\si$ are extremely mild and covered by most processes used in the literature. 

Our goal in the following is to estimate $\be$ based on irregular observations over the finite time interval $[0,1]$, say, and we will work in a setting where the observation times are typically random. In order to incorporate this additional randomness into the model we assume that the probability space contains a larger $\sigma$-field $\mathcal{G}$, and we keep using $\mathcal{F}$ to denote the $\sigma$-field with respect to which $X$ is measurable. The following condition is loosely connected with the restricted discretisation schemes introduced in Chapter 14.1 of \cite{discret} but with a slightly different predictability assumption and additional moment conditions.
\begin{condition}\label{ass:stopping}
	For each $n\in \N$ we observe the process $X$ at stopping times $0=\tau^n_0<\tau_1^n<\tau_2^n<\ldots$ with
	$\tau_0^n=0, \tau_1^n=\Delta_n\phi_1^n$ and 
	\begin{align*}
	\tau_i^n=\tau_{i-1}^n+\Delta_n \phi_{i}^n\lambda_{\tau_{i-2}^n} \text{ for all } i \ge 2
	\end{align*}
	where $\De_n \to 0$ and
	\begin{enumerate}
		\item[(a)] $\lambda_t$ is a strictly positive It\^o semimartingale w.r.t.\ the filtration $(\mathcal{F}_t)_{t\geq 0}$ and fulfills the same conditions as $\sigma_t$ stated in Assumption \ref{ass:sigma};
		\item[(b)] $(\phi_i^n)_{i\geq 1}$ is a family of random variables with respect to the $\sigma$-field $\mathcal{G}$ and independent of $\mathcal{F}$;  
		\item[(c)] $\phi_i^n\sim \phi$ for a strictly positive random variable $\phi$ with $\E[\phi]=1$, and for all $p > -2$ the moments $\E\left[\phi^p\right]$ exist.
	\end{enumerate}
\end{condition}

For all $t>0$ we define $(\fF_t^n)_{t \ge 0}$ to be the smallest filtration containing $(\fF_t)_{t \ge 0}$ and with respect to which all $\tau_i^n$ are stopping times. We also let $N_n(t)$ denote the number of observation times until $t$, i.e.\ 
	\begin{align*}
	N_n(t) = \sum_{i\geq 1}\1_{\{\tau_i^n \leq t \}},
	\end{align*}
	and of particular importance for us is the case $t=1$ because $N_n(1)$ is the (random) number of observations over the trading day $[0,1]$ from which we construct the relevant statistics later on. Note that
due to $\De_n \to 0$ we are in a high frequency situation where the time between two observations converges to zero while $N_n(1)$ diverges to infinity (both in a probabilistic sense).    

\section{Results} \label{sec:results}
\def\theequation{3.\arabic{equation}}
\setcounter{equation}{0}

The essential idea from \cite{todorov2015} is to base the estimation of the unknown activity index on the estimation of the characteristic function of a certain stable distribution. We will essentially proceed in a similar way but with some subtle changes because the underlying sampling scheme is not regular anymore. On one hand, we have to account for the fact that the time between successive observations is not constant, while on the other hand the characteristic function not only involves this particular stable distribution but also the unknown distribution $\phi$ from Assumption \ref{ass:stopping}. 

Let us become more specific here. We assume that the probability space is large enough to allow for a representation as in \cite{todorov2012}, namely that the pure-jump L\'evy process $L$ can be decomposed as
\begin{align} \label{decompL}
	L_t = S_t + \acute{S}_t- \grave{S}_t
\end{align}
where all processes on the right hand side are (possibly dependent) Lévy processes with a characteristic triplet of the form $(0,0,F)$ for a L\'evy measure of the form $F(dx) = F(x) dx$. For $S$ the L\'evy density satisfies $F(x) =  {A}|x|^{-(1+\beta)}$ while the L\'evy densities of $\acute{S}$ and $\grave{S}$ are $F(x)=|\tilde{h}(x)|$ and $F(x) = 2|\tilde{h}(x)|\1_{\{\tilde{h}(x)<0\}}$, respectively. Then $S$ is strictly $\be$-stable, 
 and its characteristic function satisfies
\begin{align*}
\E[\cos(uS_t)]=\E[\exp(iuS_t)]=\exp(-A_\beta u^\beta t)
\end{align*}
for some constant $A_\beta>0$ and any $u,t > 0.$ 

As a result of the previous decomposition (\ref{decompL}) and since $$|\tilde{h}(x)|\leq \frac{C}{|x|^{1+\beta'}}$$ for some $\be' < \be$ and all $|x|\leq x_0$ holds due to Condition \ref{ass:general}, it is clear that the jump behaviour of $L$ and thus of $X$ is governed by the $\be$-stable process $S$ for small time intervals. This observation is the key to our following estimation procedure: Based on the high-frequency observations of $X$ we will first estimate a function $L(p,u,\beta)$ which, as noted before, is related to the characteristic function of $S$ but involves the unknown distribution of $\phi$ as well. Here, $p$ and $u$ are additional parameters that can be chosen by the statistician. In a second step we will essentially use a Taylor expansion of $L$ (as a function of $u$) around zero to finally come up with an estimator for $\beta$. 

In the following, we denote with 
$$\widetilde{\Delta_i^nX}=\frac{\Delta_n}{\tau_i^n-\tau_{i-1}^{n}}\left(X_{\tau_i^n}-X_{\tau_{i-1}^{n}}\right)$$
the $i$th increment of the process $X$, but where we have included an additional rescaling in order to account for the different lengths of the intervals in an irregular sampling scheme, and we occasionally also use $\widetilde{\Delta_i^n S}$ for the rescaled increment of the $\be$-stable $S$. For any $p > 0$ and $u > 0$ we then set
\begin{align*}
\widetilde{L}^n(p,u)=\frac{1}{\Nn-k_n-2}\sum_{i=k_n+3}^{\Nn}\cos\left(u\frac{\widetilde{\Delta_i^n X}-\widetilde{\Delta_{i-1}^nX} }{(\widetilde{V}_i^n(p))^{1/p}}\right)
\end{align*}
where the auxiliary sequence $k_n$ satisfies $k_n \to \infty$ and $k_n \De_n \to 0$ and where 
\begin{align*}
\widetilde{V}_i^n(p)=\frac{1}{k_n}\sum_{j=i-k_n-1}^{i-2}|\widetilde{\Delta_j^n X}-\widetilde{\Delta_{j-1}^n X}|^p, ~~i=k_n+3,\ldots,\Nn,
\end{align*}
is used to estimate the unknown local volatility $\si$. 

At first, it seems somewhat odd to include $\De_n$ in the definition of $\widetilde{\Delta_i^nX}$ because this quantity cannot be observed in practice. We will base our statistical procedure in the following on $\widetilde{L}^n(p,u)$, however, and it is obvious from its definition that it is in fact independent of $\De_n$ as the latter appears as a factor both in the numerator and in the denominator. Thus we are safe to work with $\widetilde{\Delta_i^nX}$, and its definition makes it easier to compare its results with the standard increment $\De_i^n X = X_{\tau_i^n}-X_{\tau_{i-1}^{n}}$. These obviously coincide in the case of a regular sampling scheme. Note also that, even though its asymptotic condition is stated in terms of $\De_n$, the choice of $k_n$ can in practice be based on the size of $N_n(1)$ which essentially grows as $\De_n^{-1}$.

The main part of the upcoming analysis is devoted to the study of the asymptotic behaviour of $\widetilde{L}^n(p,u)$. Its definition together with the previous discussion suggests that its limit should involve the characteristic function of $S$, but it also becomes apparent that the limit cannot be independent of $\phi$. We will prove in the following that the first order limit is 
\begin{align*}
L(p,u,\beta)=\E \left[\exp \left(-u^\beta C_{p,\beta}((\phi^{(1)})^{1-\beta}+(\phi^{(2)})^{1-\beta})\right) \right]
\end{align*} 
with the constant $C_{p,\beta}$ being defined via 
\begin{align} \label{defcon}
\mu_{p,\beta}=\E[|S_1|^p]^{\frac{\beta}{p}},\quad 	\kappa_{p,\beta}=\E\left[\left((\phi^{(1)})^{1-\beta}+(\phi^{(2)})^{1-\beta}\right)^\frac{p}{\beta}\right]^\frac{\beta}{p}, \quad C_{p,\beta}=\frac{A_\beta}{\mu_{p,\beta}\kappa_{p,\beta}}>0,
\end{align}
and where $\phi^{(1)}$ and $\phi^{(2)}$ denote two independent copies with the same distribution as $\phi$, defined on an appropriate probability space. For simplicity, we still use $\E[\cdot]$ to denote the expectation on this generic space. The first main theorem then reads as follows: 

\begin{theorem} \label{thmlncon}
Suppose that Conditions \ref{ass:general}--\ref{ass:stopping} are in place and let $k_n\sim C_1 \De_n^{-\varpi}$ for some $C_1> 0$ and some $\varpi\in(0,1)$. Then we have 
\[
\widetilde{L}^n(p,u) \pn L(p,u,\beta)
\]
for any fixed $u > 0$ and any choice of $0 < p < \be/2$. 		
\end{theorem}

While this result is interesting in itself, at first glance it does not help much for the estimation of $\beta$ because $L(p,u,\beta)$ depends in a complicated way on the unknown distribution of $\phi$. If we utilize the familiar approximation $\exp(y) = 1+y+ o(y)$ for $y \to 0$, however, it seems reasonable to hope that the approximation
\begin{align} \label{approxL}
L(p,u,\beta) \approx 1 + \E\left[-u^\beta C_{p,\beta}((\phi^{(1)})^{1-\beta}+(\phi^{(2)})^{1-\beta})\right] = 1 - u^\beta C_{p,\beta}\kappa_{\beta,\beta}
\end{align}
holds for any choice of a small $u > 0$, which now is much easier to handle. Namely, an estimator for $\beta$ is then based on an appropriate combination of two estimators $\widetilde{L}^n(p,u_n)$ and $\widetilde{L}^n(p,v_n)$ with $u_n \to 0$ and $v_n = \rho u_n$ for some $\rho > 0$. Precisely, we set 
\begin{align*}
\hat{\beta}(p,u_n,v_n)=\frac{\log(-(\widetilde{L}^n(p,u_n)-1))-\log(-(\widetilde{L}^n(p,v_n)-1))}{\log(u_n/v_n)}
\end{align*}
which obviously is symmetric upon exchanging $u_n$ and $v_n$. 

\begin{remark} \label{be2}
We will often choose $\rho = 1/2$ in which case $\hat{\beta}(p,u_n,v_n) \le 2$ can be shown. This is of course a desirable property as it resembles the bound for the stability index $\beta$ itself, but it also bears some restrictions regarding the quality of a limiting normal approximation for values of $\beta$ close to 2. See Figures \ref{fig:beta11} and \ref{fig:beta12} below.
\end{remark}

Before we discuss the asymptotic behaviour of the estimator $\hat{\beta}(p,u_n,v_n)$ we will state a bivariate central limit theorem for $\widetilde{L}^n(p,u_n) - L(p,u_n,\beta)$ and $\widetilde{L}^n(p,v_n) - L(p,v_n,\beta)$ with $u_n$ and $v_n$ chosen as above.  

\begin{theorem}\label{thmlnclt}
Suppose that Conditions \ref{ass:general}--\ref{ass:stopping} are in place and let $k_n\sim C_1 \De_n^{-\varpi}$ for some $C_1> 0$ and some $\varpi\in(0,1)$ as well as $u_n\sim C_2 \De_n^{\varrho}$ for some $C_2 > 0$ and $\varrho \in (0,1)$. Suppose further that
\begin{align*}
	\beta'<\frac{\beta}{2}, \quad \frac{1}{3}\vee\frac{1}{8\varrho}<p<\frac{\beta}{2}, \quad \varpi\geq\frac{2}{3}, \quad \frac{1}{3\beta}<\varrho< \frac{1}{\beta}, \quad \frac{1}{\beta}<\frac{\varpi}{p}-\varrho, \quad 2\varpi-\varrho\beta<1 
	\end{align*}
	hold. Then 
\begin{align*}
\left(\frac{\sqrt{N_n(1)}}{u_n^{\beta/2}}(\widetilde{L}^n(p,u_n)-L(p,u_n,\beta)),\frac{\sqrt{N_n(1)}}{v_n^{\beta/2}}(\widetilde{L}^n(p,v_n)-L(p,v_n,\beta))\right)
\end{align*}
converges $\fF$-stably in law to a limit $(X',Y')$ which is jointly normal distributed (independent of $\fF$) with mean $0$ and covariance matrix $\mathcal{C'}$ given by
	\begin{align*}
	\mathcal{C}_{11}'=\mathcal{C}_{22}'=C_{p,\beta}\kappa_{\beta,\beta}(4-2^\beta), \quad \mathcal{C}_{12}'=\mathcal{C}_{21}'=C_{p,\beta}\kappa_{\beta,\beta}\frac{2+2\rho^\beta - (1+\rho)^\beta  -|1-\rho|^\beta }{\rho^{\beta/2}}.
	\end{align*}
\end{theorem}

\begin{remark} \label{rem:feas}
	The above choice of the parameters $\varrho$, $\varpi$ and $p$ is feasible even if we do not know $\beta$. It can easily be seen that e.g.\ $\varrho=\frac{1}{3}, \varpi=\frac23$ and any $p \in (\frac38,\frac12)$ satisfies the conditions in Theorem \ref{thmlnclt}. 
\end{remark}

The following result is the main theorem of this work and it provides the central limit theorem for the estimator $\hat{\beta}(p,u_n,v_n)$. Its proof builds heavily on Theorem \ref{thmlnclt}. 

	\begin{theorem}\label{thm:beta_hat}
	Under the conditions of Theorem \ref{thmlnclt} we have the $\fF$-stable convergence in law
	\begin{align*}
	u_n^{\beta/2}\sqrt{N_n(1)}(\hat{\beta}(p,u_n,v_n)-\beta)\tols X
	\end{align*} 
	where $X$ is a normal distributed random variable (independent of $\fF$) with mean $0$ and variance $$\frac{(\rho^\beta+1)(4-2^\beta)-2(2+2\rho^\beta - (1+\rho)^\beta  -|1-\rho|^\beta )}{\kappa_{\beta,\beta}\rho^\beta \log(1/\rho)^2 C_{p,\beta}}.$$
\end{theorem}

A simple corollary is the consistency of $\hat{\beta}(p,u_n,v_n)$ as an estimator for $\beta$. 

	\begin{corollary}\label{cor:beta_hat}
	Under the conditions of Theorem \ref{thmlnclt} we have
	\begin{align*}
	\hat{\beta}(p,u_n,v_n) \pn \beta.
	\end{align*} 
\end{corollary}

For a feasible application of Theorem \ref{thm:beta_hat} we need a consistent estimator for the variance of the limiting normal distribution which essentially boils down to the estimation of $\ka_{\be, \be}$. This problem will be discussed in the next section, alongside with a thorough analysis of the finite sample properties of $\hat{\beta}(p,u_n,v_n)$.

\section{Simulation study} \label{sec:simul}
\def\theequation{4.\arabic{equation}}
\setcounter{equation}{0}

\begin{table} 
\begin{center}

		\begin{tabular}{|c|c|c|c|c|}
			\hline 
			$\beta$	& Mean of $\hat{\beta}(p,u,v)$    & Empirical Variance  &  Theoretical Variance \\ 
			\hline 
			$1.1$	& $1.1181~(1.0455)$   & $7.2689~(3.2672)$ & $7.2457~(3.3802)$ \\ 
			\hline 
			$1.3$	& $1.3123~(1.2356)$   & $5.2131~(2.2777)$ & $5.4853~(2.2274)$ \\ 
			\hline 
			$1.5$	& $1.4923~(1.4421)$ & $3.2153~(1.3201)$  & $3.907~(1.3817)$  \\ 
			\hline 
			$1.7$	& $1.7173~(1.6086)$   & $1.6501~(0.73274)$  & $2.354~(0.7245)$ \\ 
			\hline 
			$1.9$	& $1.8849~(1.7759)$   & $0.425~(0.2716)$ & $0.7852~(0.2107)$ \\ 
			\hline 
		\end{tabular} 	
		\caption{Results for $\rho=1/2$ and $\De_n^{-1}=1000$ are shown. The second column shows the empirical mean of the $N=1000$ samples of $\hat{\beta}(p,u_n,v_n)$. The third one is the empirical variance of $u_n^{\beta/2}\sqrt{N_n(1)}(\hat{\beta}(p,u_n,v_n)-\beta)$, whereas the final column shows the true asymptotic variance from Theorem \ref{thm:beta_hat}. In brackets the same results are given for $\rho=2$.} \label{tab1}
\end{center}
\end{table}

This chapter deals with the numerical assessment of the finite sample properties of $\hat{\beta}(p,u_n,v_n)$, and we also include a discussion regarding the estimation of the variance in the central limit theorem in order to obtain a feasible result. In the following, let $W$ be a standard Brownian motion and $L$ a symmetric stable process with a Lévy density
\begin{align*}
h(x)=\frac{1}{|x|^{1+\beta}}
\end{align*}
for some $\beta\in(1,2)$. We then set 
\begin{align*}
	\alpha_t = \int_{0}^{t} 2(1-\alpha_s) ds + 2 \int_{0}^{t} d{W}_s, \quad \sigma_t = \int_{0}^{t} \alpha_s d{W}_s, 
\end{align*}
and assume that we observe
\begin{align*}
X_t=X_0+\int_{0}^{t}\alpha_s ds+\int_{0}^{t}\sigma_{s-}dL_s
\end{align*}
which obviously fulfills Conditions \ref{ass:general} and \ref{ass:sigma}. For the observation scheme we choose
\begin{align}
	\lambda_t = \int_{0}^{t} (5-\lambda_s) ds + \int_{0}^{t} d\widetilde{W}_s, \quad \phi = \frac{\phi' \vee 0.1}{\E\left[\phi' \vee 0.1\right]} \label{num:phi},	
\end{align}
with $\phi' \sim \text{Exp}(1)$ and with the starting values of the processes being $\alpha_0 = \sigma_0 = X_0 = \lambda_0 = 1$. We assume that $\widetilde{W}$ is a standard Brownian motion as well, independent of $W$. The purpose of the minimum in the definition of $\phi$ in \eqref{num:phi} is to ensure the (negative) moment condition from Assumption \ref{ass:stopping} (c) to hold, which (as can be seen from additional simulations) not only seems to be relevant in theory but in practice as well. Note also that the choice of $\lambda_0=1$ combined with the mean reversion of $\lambda$ to $5$ leads to pronounced changes in the distribution of the $\tau_i^n$ over time. 

\begin{table} 
\begin{center}
		\begin{tabular}{|c|c|c|c|c|}
			\hline 
			$\beta$	& Mean of $\hat{\beta}(p,u,v)$    & Empirical Variance  &  Theoretical Variance \\ 
			\hline 
			$1.1$	& $1.1048~(1.0852)$   & $7.451~(3.3392)$ & $7.2458~(3.3802)$ \\
			\hline 
			$1.3$	& $1.3067~(1.2763)$   & $5.257~(2.229)$ & $5.4853~(2.2274)$ \\ 
			\hline 
			$1.5$	& $1.5197~(1.4772)$ & $3.6036~(1.3626)$  & $3.907~(1.3817)$  \\ 
			\hline 
			$1.7$	& $1.7088~(1.6743)$  & $2.1925~(0.7034)$  & $2.354~(0.7245)$ \\ 
			\hline 
			$1.9$	& $1.914~(1.8726)$   & $0.4864~(0.21937)$ & $0.7852~(0.2107)$ \\ 
			\hline 
		\end{tabular} 
		\caption{Results for $\rho=1/2$ and $\De_n^{-1}=10,000$ are shown. The second column shows the empirical mean of the $N=1000$ samples of $\hat{\beta}(p,u_n,v_n)$. The third one is the empirical variance of $u_n^{\beta/2}\sqrt{N_n(1)}(\hat{\beta}(p,u_n,v_n)-\beta)$, whereas the final column shows the true asymptotic variance from Theorem \ref{thm:beta_hat}. In brackets the same results are given for $\rho=2$.} \label{tab2}
\end{center}
\end{table}

For the simulation of $X$ we use a standard Euler scheme, utilising Proposition 1.7.1 in \cite{stableprocesses1994} to obtain symmetric stable random variables. We also set $u_n = \Nn^{-1/3}$, $k_n = \Nn^{2/3}$, $p=1/2$, essentially in accordance with Remark \ref{rem:feas}. Below, we present the results for $\beta\in\{1.1,1.3,1.5,1.7,1.9\}$ and $\rho\in\{1/2,2\}$, generating $N = 1000$ samples, and we discuss both $\De_n^{-1} = 1000$ and $\De_n^{-1} = 10,000$. Note that our choice of $\la$ in \eqref{num:phi} yields about $\Nn\approx520$ observations in the first case, whereas $\Nn\approx5200$ in the second one.

Table \ref{tab1} shows mixed results but nevertheless allows to draw some conclusions: First, we see for a choice of $\rho=1/2$ that the estimator for $\be$ behaves correctly on average whereas the (relative) difference between the true and the estimated variances grows with $\be$. Note that the larger choice of $\rho$ effectively corresponds to choosing $u_n$ twice as large. Hence Table \ref{tab1} confirms empirically what is already known from the construction of the estimator: It relies on $u_n \to 0$, so a larger choice of $\rho$ induces an additional bias. On the other hand, the rate of convergence to the normal distribution improves as $u_n$ becomes larger, and it also seems as if the quality of the variance estimation is better in this case. 

\begin{center}
\begin{figure}
\centering
\includegraphics[width=0.4\linewidth]{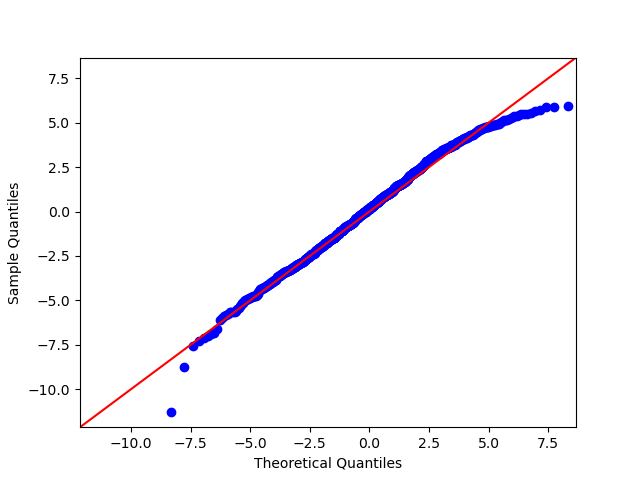}\includegraphics[width=0.4\linewidth]{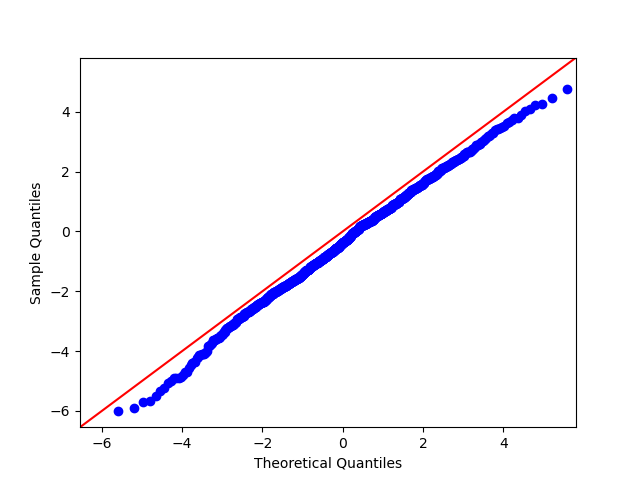}
\includegraphics[width=0.4\linewidth]{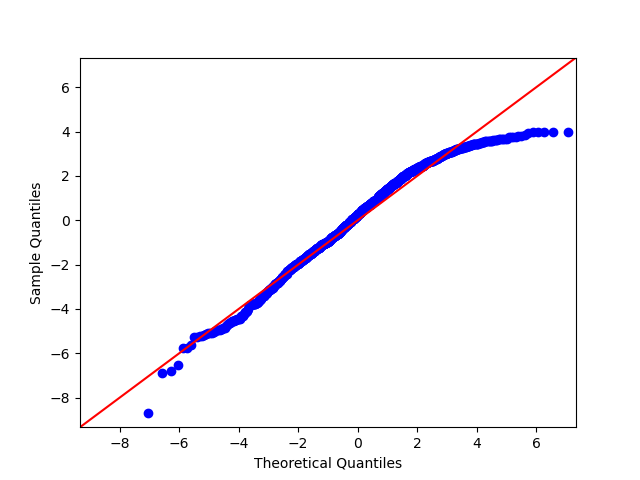}\includegraphics[width=0.4\linewidth]{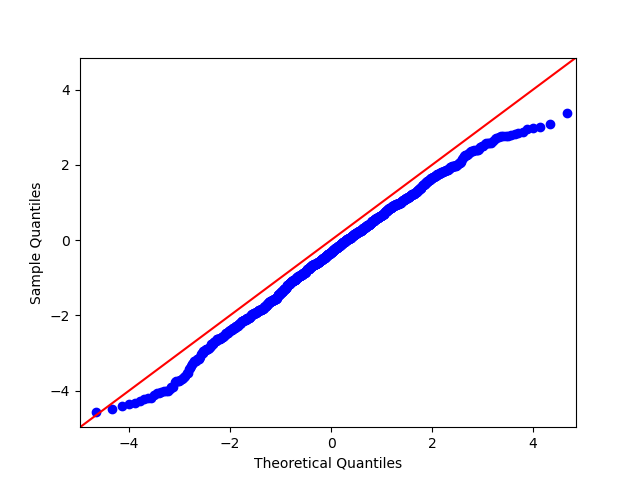}
\includegraphics[width=0.4\linewidth]{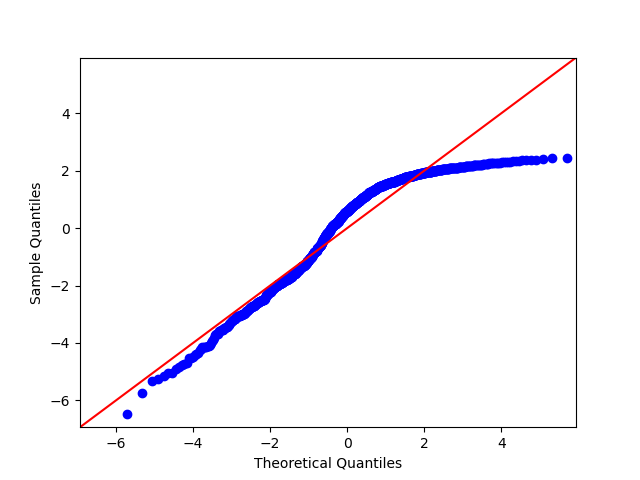}\includegraphics[width=0.4\linewidth]{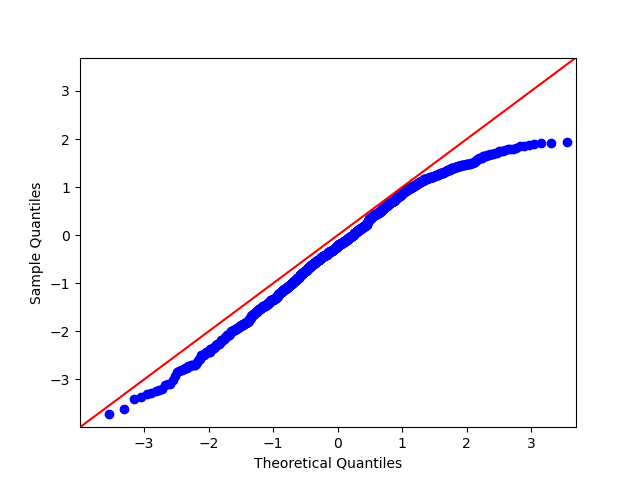}
\includegraphics[width=0.4\linewidth]{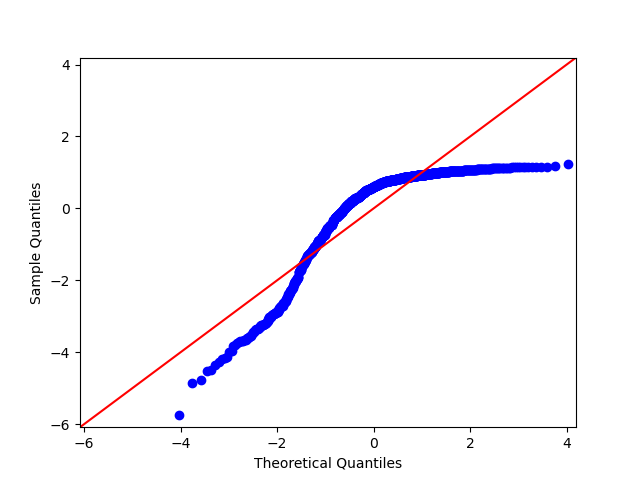}\includegraphics[width=0.4\linewidth]{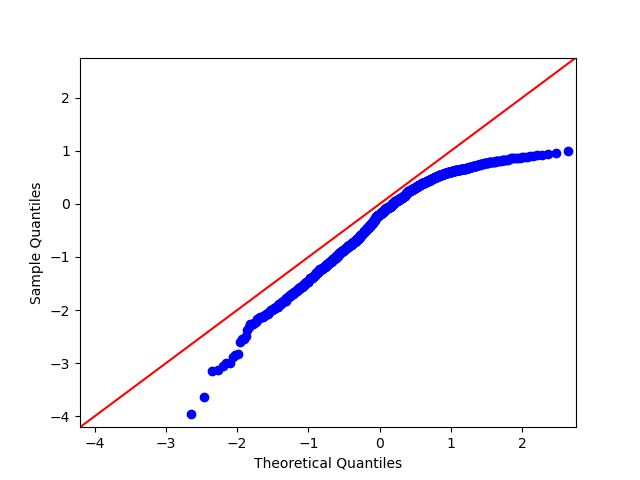}
\includegraphics[width=0.4\linewidth]{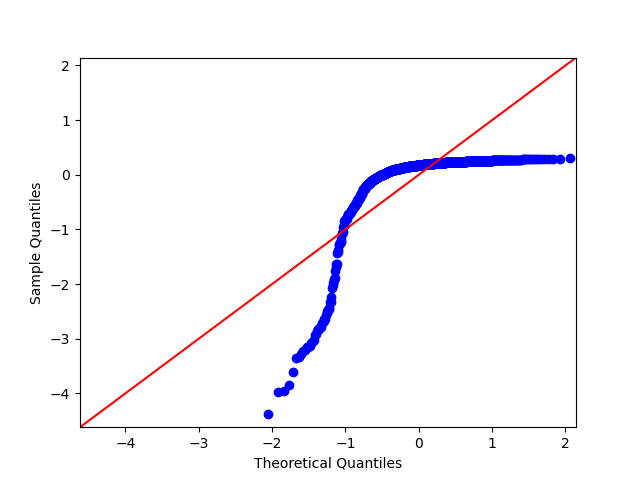}\includegraphics[width=0.4\linewidth]{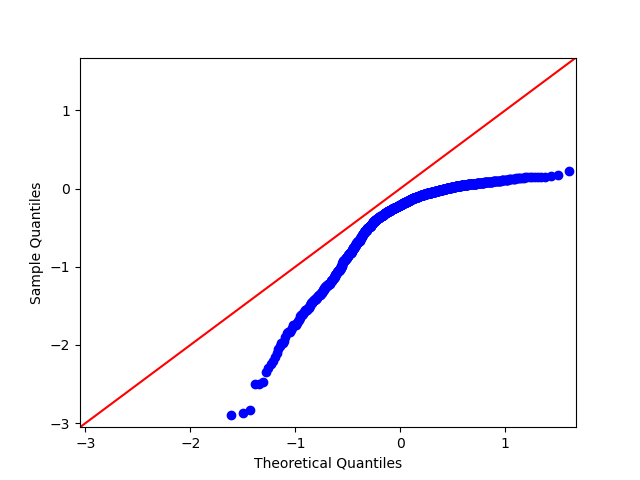}
\caption{QQ-plots in the case $N=1000$, $\De_n^{-1}=1000$ and $\beta \in \{1.1,1.3,1.5,1.7,1.9\}$ are shown, with $\beta$ increasing from the top to the bottom. The left hand side provides the plots for $\rho=0.5$, the right hand side the ones for $\rho=2$.}
\label{fig:beta11}
\end{figure}
\end{center}
\begin{figure}
	\centering
	\includegraphics[width=0.4\linewidth]{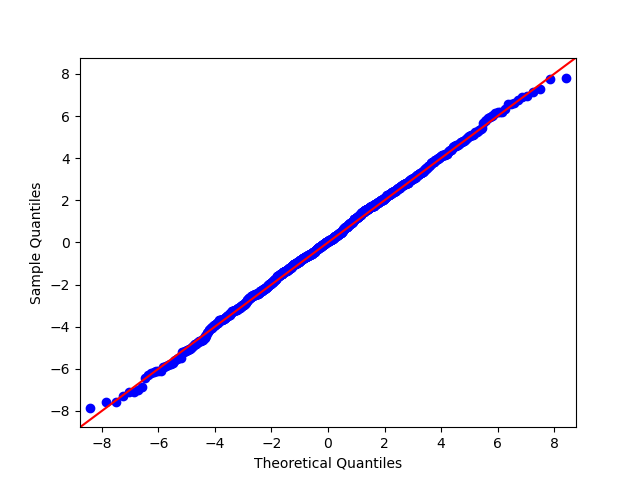}\includegraphics[width=0.4\linewidth]{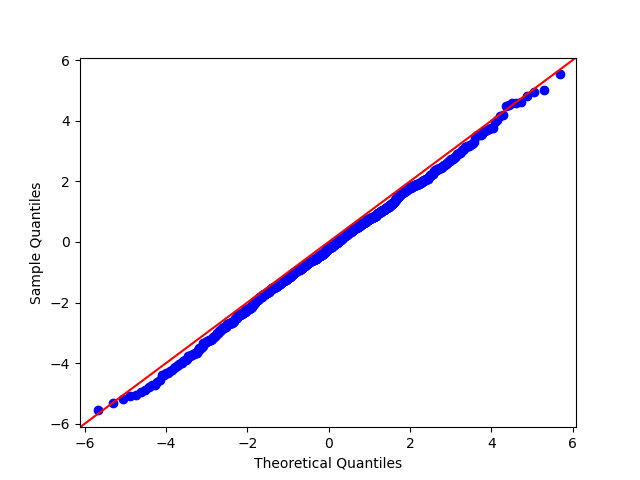}
	\includegraphics[width=0.4\linewidth]{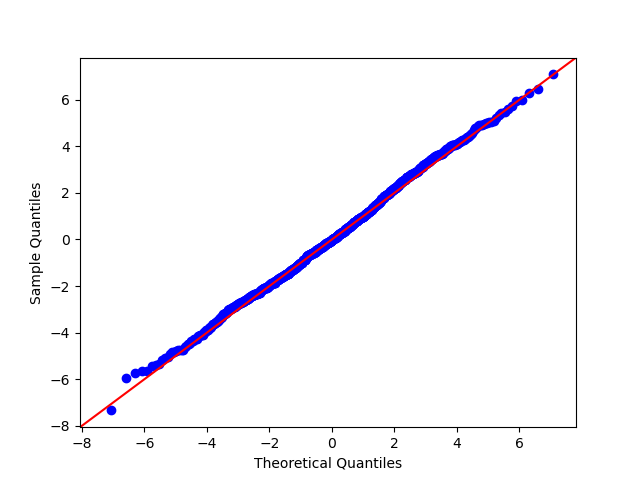}\includegraphics[width=0.4\linewidth]{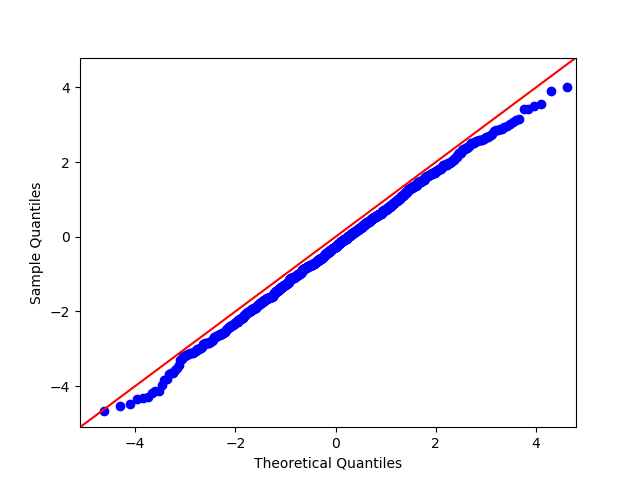}
	\includegraphics[width=0.4\linewidth]{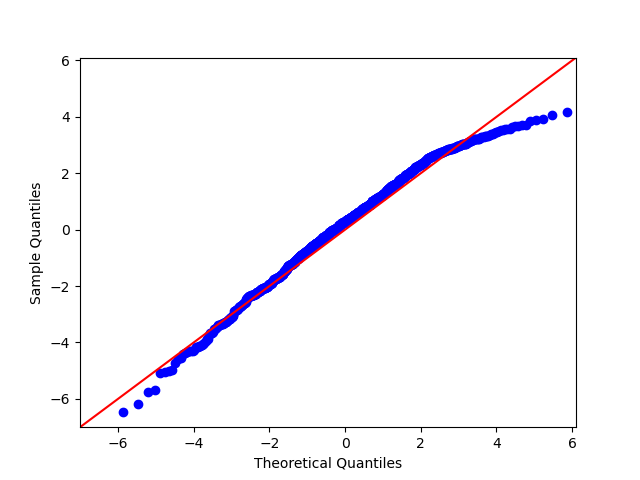}\includegraphics[width=0.4\linewidth]{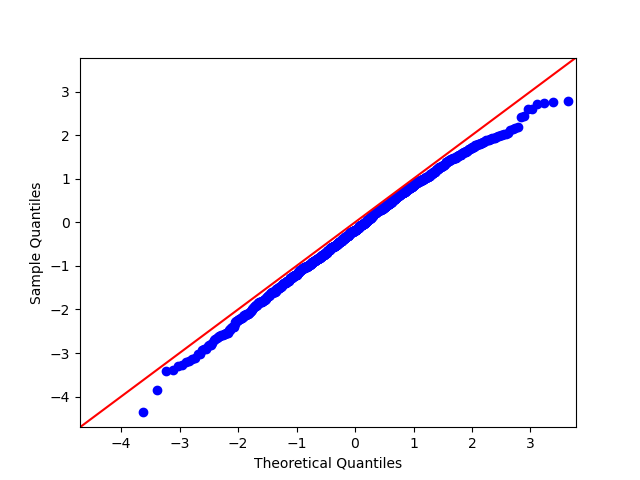}
	\includegraphics[width=0.4\linewidth]{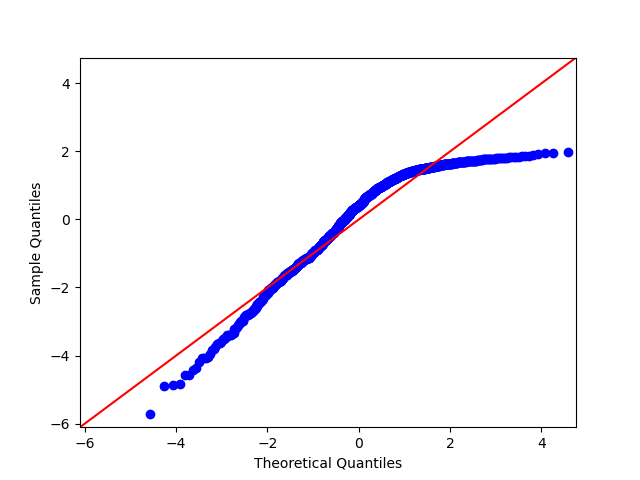}\includegraphics[width=0.4\linewidth]{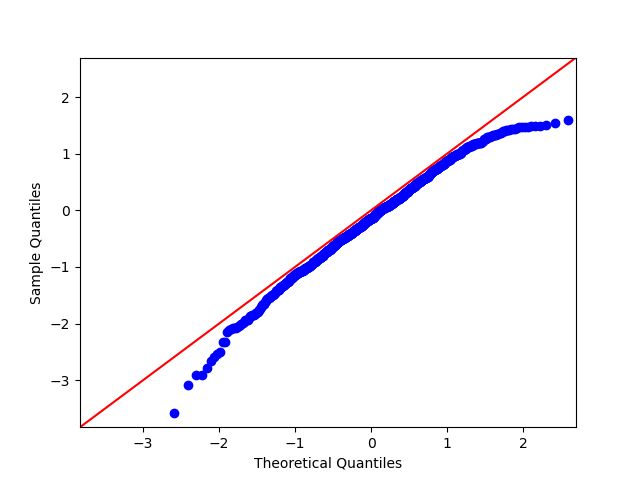}
	\includegraphics[width=0.4\linewidth]{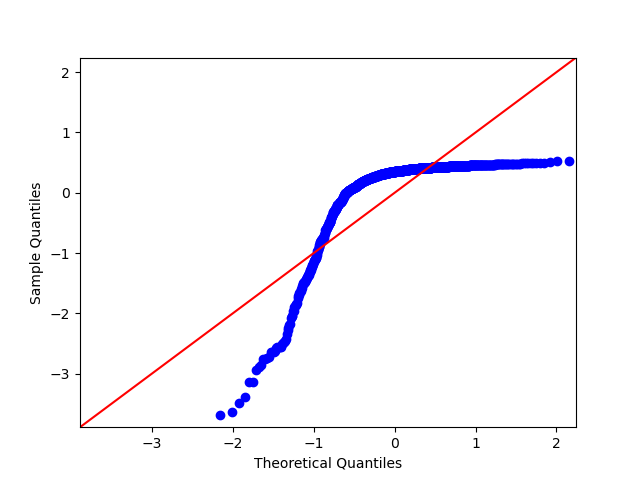}\includegraphics[width=0.4\linewidth]{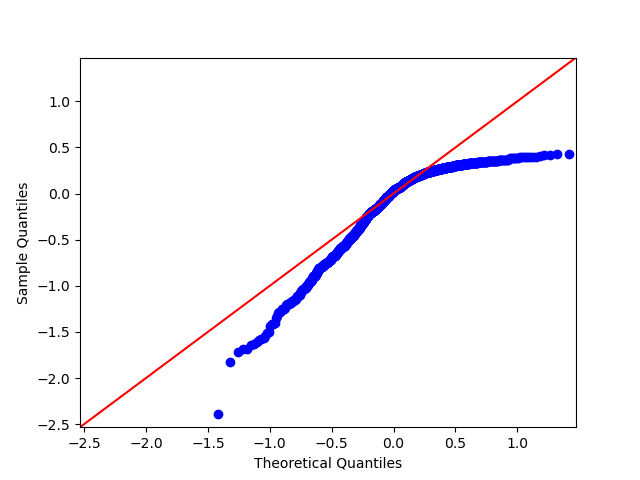}
\caption{QQ-plots in the case $N=1000$, $\De_n^{-1}=10,000$ and $\beta \in \{1.1,1.3,1.5,1.7,1.9\}$ are shown, with $\beta$ increasing from the top to the bottom. The left hand side provides the plots for $\rho=0.5$, the right hand side the ones for $\rho=2$.}
	\label{fig:beta12}
\end{figure}

We follow this discussion with Table \ref{tab2} which is constructed in the same manner as above for $\De_n^{-1}=10,000$, and we basically see improvement across the board. Now, the bias for both $\rho=1/2$ and $\rho=2$ is very small, even for large values for $\beta$, and we can also observe that the approximated variance specifically for $\beta\in\{1.5,1.7,1.9\}$ is much closer to the theoretical one than previously.

In a second step we present some QQ-plots to visualise the quality of the approximating normal distribution from Theorem \ref{thm:beta_hat}. In this case we need to estimate the limiting variance 
\[
\frac{(\rho^\beta+1)(4-2^\beta)-2(2+2\rho^\beta - (1+\rho)^\beta  -|1-\rho|^\beta )}{\kappa_{\beta,\beta}\rho^\beta \log(1/\rho)^2 C_{p,\beta}}
\]
in Theorem \ref{thm:beta_hat}, and we note from eq.\ (3.3) in \cite{todorov2015} that 
	\begin{align*}
	\frac{A_\beta}{\mu_{p,\beta}} = \left(\frac{2^p\Gamma((1+p)/2)\Gamma(1-p/\beta)}{\sqrt{\pi}\Gamma(1-p/2)}\right)^{-\beta/p}
	\end{align*}
holds. Hence, the only unknown quantities in the variance are $\be$ and $\kappa_{\be, \be}$, and from Corollary \ref{cor:beta_hat} we only need a consistent estimator for $\ka_{\be, \be}$ to obtain a consistent plug-in estimator for the limiting variance. 

By definition, $\kappa_{\be,\beta}=\E\left[\left((\phi^{(1)})^{1-\beta}+(\phi^{(2)})^{1-\beta}\right)\right]$. Hence, using $\tau_i^n-\tau_{i-1}^n = \Delta_n \phi_{i}^n\lambda_{\tau_{i-2}^n}$ and the c\`adl\`ag property of $\la$, a natural way to construct an estimator for $\kappa_{\be, \be}$ is to build it from sums of adjacent increments of the $\tau_i^n$, rescaled by the length of the total time interval relative to the number of increments. Whenever $\be$ needs to be included, it is replaced by its consistent estimator $\hat{\beta}_n = \hat \beta(p, u_n, v_n)$. The consistency of such an estimator for $\ka_{\be, \be}$ is formally given in the following lemma. Its proof is rather straightforward but lengthy and therefore omitted. 
\begin{lemma} \label{lemvarest}
	Let $r_n \sim C_3 \De_n^{-\psi}$ for some $\psi \in (0,1)$. Then
	\begin{align*}
	\widehat{\kappa}_n=\frac{1}{\Nn-r_n-2}\sum_{i=r_n+3}^{\Nn} \chi_i^n \pn \kappa_{\beta,\beta}
	\end{align*}
	where 
	\begin{align*}
	\chi_i^n = \left(\left(\frac{r_n}{\tau_{i-2}^n-\tau_{i-2-r_n}^n}\right)^{1-\hat{\beta}_n}\left(\left(\tau_{i}^n-\tau_{i-1}^n\right)^{1-\hat{\beta}_n}+\left(\tau_{i-1}^n-\tau_{i-2}^n\right)^{1-\hat{\beta}_n}\right)\right).
	\end{align*}
\end{lemma}

We use the same configuration of parameters as discussed earlier and further set $r_n = \Nn^{4/5}$. Due to the choice of $\rho=1/2$ and $\rho=2$ Remark \ref{be2} applies. As noted before, this condition prevents the normal approximation from working well in the right tails, and as expected Figure \ref{fig:beta11} shows that this effect becomes more pronounced as $\beta$ gets closer to the upper bound. For $\De_n^{-1}=1000$, this dubious tail behaviour already starts to appear for $\beta=1.5$. Nevertheless, it should be noted that the quality of the distributional approximation increases visibly with the higher sample-size $\De_n^{-1}=10,000$ for both choices of $\rho$. Figure \ref{fig:beta12} shows for instance that $\be=1.5$ is not really critical anymore, i.e.\ an increasing sample size allows for an accurate approximation of larger values of $\be$. Also, a slight improvement from $\rho=1/2$ to $\rho=2$ can be noted, in line with the previous discussion regarding the rate of convergence. 

A natural way to further improve the finite sample properties is to conduct a bias correction for the higher order terms. Our current approach to estimate $\beta$ relies on the approximation (\ref{approxL}) while a more precise one would e.g.\ be a third order expansion of the form
\begin{align*}
& L(p,u_n,\be) -1 \\ =& C_{p,\beta} u_n^\beta \left(-\kappa_{\beta,\beta} +\E\left[ \frac{C_{p,\beta}u_n^\beta ((\phi^{(1)})^{1-\beta}+(\phi^{(2)})^{1-\beta})^2}{2}\right] - \E\left[ \frac{(C_{p,\beta}u_n^\beta)^2 ((\phi^{(1)})^{1-\beta}+(\phi^{(2)})^{1-\beta})^3}{6}\right]  \right).
\end{align*}
An estimator for $\be$ is then given by
\begin{align*}
\overline{\beta}(p,u_n,v_n)=\frac{\log(-(\widetilde{L}^n(p,u_n)-1 - \widehat{deb}(\beta, r, u_n, C_{p,\beta})))-\log(-(\widetilde{L}^n(p,v_n)-1 - \widehat{deb}(\beta, r, v_n, C_{p,\beta})))}{\log(u_n/v_n)}
\end{align*}
where $\widehat{deb}(\beta, r, u_n, C_{p,\beta})$ estimates
\[
\E\left[ \frac{(C_{p,\beta}u_n^\beta ((\phi^{(1)})^{1-\beta}+(\phi^{(2)})^{1-\beta}))^2}{2}\right] - \E\left[ \frac{(C_{p,\beta}u_n^\beta ((\phi^{(1)})^{1-\beta}+(\phi^{(2)})^{1-\beta}))^3}{6}\right]
\]
and can be constructed analogously to Lemma \ref{lemvarest}. Ideally, such a correction would allow for a bigger choice of $u_n$ in finite samples, thus leading to a better rate of convergence. 

\begin{figure}
	\centering
	\includegraphics[width=0.4\linewidth]{rho12n1kbeta=1,7.png}\includegraphics[width=0.4\linewidth]{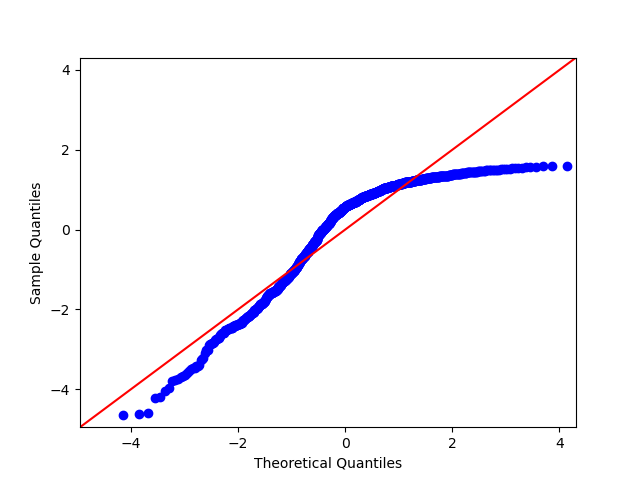}
\caption{QQ-plots in the case $N=1000$, $\De_n^{-1}=1000$, $\beta =1.7$ and $\rho = 0.5$. On the left: The standard estimator $\hat{\beta}(p,u_n,v_n)$. On the right: The bias-corrected estimator $\overline{\beta}(p,u_n,v_n)$.}
	\label{fig:beta13}
\end{figure}

As an example, we discuss the case of $\De_n^{-1}=1000$ with $\be = 1.7$ and $\rho = 0.5$, but we now chooce $u_n = \Nn^{-0.28}$ and keep all other variables unchanged. In this case, mean and empirical variance become $1.7042$ and $1.7970$, both improving the corresponding values from Table \ref{tab1}. Also, the corresponding QQ-plot in Figure \ref{fig:beta13} clearly shows a better approximation of the limiting normal distribution, still with the original problems in the right tail. 

\section{Proofs} \label{sec:proof}
\def\theequation{5.\arabic{equation}}
\setcounter{equation}{0}

\subsection{Prerequisites on localisation}

As usual one starts with localisation results, i.e.\ with results that allow to prove the main theorems under conditions which are slightly stronger than Conditions \ref{ass:general} and \ref{ass:sigma} for the processes involved and also stronger than Condition \ref{ass:stopping} on the sampling scheme. We begin with the additional assumptions on the processes. 

\begin{condition}\label{ass:strong}
In addition to Conditions \ref{ass:general} and \ref{ass:sigma} we assume that
	\begin{itemize}
		\item[(a)] $|\sigma_t|$ and $|\sigma_t|^{-1}$ are uniformly bounded;
		\item[(b)] $|\delta^\alpha(t,x)|+|\delta^\sigma(t,x)|\leq \gamma(x)$ for all $t>0$, where $\gamma(x)$ is a deterministic bounded function on $\R$ with $\int_E|\gamma(x)|^r\lambda(dx)<\infty$ for some $0\leq r<2$;
		\item[(c)] $b^\alpha, b^\sigma, \eta^\alpha, \eta^\sigma, \widetilde{\eta}^\alpha, \widetilde{\eta}^\si, \overline{\eta}^\alpha$ and $\overline{\eta}^\si$ are bounded;
		\item[(d)] the process $\left(\int_{\R}\left(|x|^{\beta'}\wedge 1\right)\nu_t^Y(dx)\right)_{t\geq 0}$ is bounded and the jumps of $Y$ are bounded;
		\item[(e)] the jumps of $\acute{S}$ and $\grave{S}$ are bounded. 
	\end{itemize}
The same properties hold for the processes that govern $\la$. 	
\end{condition} 

The following lemma gives the formal result why we can assume in the following that the strengthened Condition \ref{ass:strong} holds, namely because we are interested in $X$ on the bounded interval $[0,1]$ only and eventually $E_p > 1$ for a localising sequence, at least with a probability converging to 1. Its proof closely resembles the one of Lemma 4.4.9 in \cite{discret} which is why we refer the reader to part 3) of their proof. 
	
\begin{lemma}\label{lem:stopping}
	Let $X$ be a process fulfilling Condition \ref{ass:general} and \ref{ass:sigma}. Then, for each $p>0$ there exists a stopping time $E_p$ and a process $X(p)$ such that $X(p)$ and its components, $\alpha(p)$, $\sigma{(p)}$ and $Y(p)$, fulfill Assumption \ref{ass:strong}, and it also holds that  $X(p)_t = X_t$ for all $t<E_p$. The sequence of stopping times can be chosen such that $E_p\nearrow\infty$ almost surely when $p\rightarrow\infty$.
\end{lemma}

For all proofs concerning the asymptotics of $\widetilde{L}^n(p,u)$ and $\hat{\beta}(p,u_n,v_n)$ it becomes important that the process $\lambda_t$ driving the observation times $\tau_i^n$ is bounded from above and below. This means that we need a stronger assumption than just Condition \ref{ass:stopping} as well, and we also need to assume that for a given $n$ the number of observations until any fixed $T$ is bounded by a constant times $\De_n^{-1} T$. 
\begin{condition}\label{ass:strongstopping}
	In addition to Condition \ref{ass:stopping} there exists some $C>1$ such that 
	\begin{itemize}
		\item[(a)] The process $\lambda$ fulfills the same assumptions as $\sigma$ in Condition \ref{ass:strong}, and in particular we have for all $t>0$ 
		\begin{align*}
		\frac{1}{C} \leq \lambda_t \leq C.
		\end{align*}
		\item[(b)] For any given $n$ and any $T > 0$ we have
		\begin{align*}
		N_n(T) \leq C \De_n^{-1} T.
		\end{align*}
	\end{itemize}
\end{condition}

Strengthening Condition \ref{ass:stopping} ultimatively results in changing the entire observation scheme which makes it somewhat harder to formally prove that such an assumption is indeed adequate. We begin with a result on the boundedness of $\la$ as in part (a) above, and for every $n$ let $F_n$ be a random variable which not just depends on $n$ but also on the process $X$ and on the discretisation scheme via $\lambda$ and the variables $\phi_i^n$. Likewise, a possible stable limit $F$ of $F_n$ is assumed to depend on the same factors and is realised on an extension $(\widetilde{\Omega},\widetilde{\mathcal{G}},\widetilde{\mathbb{P}})$ of the original probability space $(\Omega,\mathcal{G},\mathbb{P})$. Furthermore, for each $C > 1$ we define $\lambda^{(C)}_t$ in such a way that $\frac{1}{C} \leq \lambda^{(C)}_t \leq C$ holds. We then set $E_C$ to be the stopping time from Lemma \ref{lem:stopping} with $C$ replacing $p$, and this lemma can be applied because, by Condition \ref{ass:stopping}, $\la$ is assumed to satisfy the same structural properties as $\si$.

\begin{lemma}\label{Lem:Stop1}
	Assume that Assumption \ref{ass:stopping} holds and construct, for each $C>1$, each stopping time $E_C$ and each process $\lambda^{(C)}$, a new discretisation scheme, i.e.\ new stopping times $\{\tau_i^{n,C}:i \geq 0\}$ and a new $N_n^C(T)$ as in Condition \ref{ass:stopping} but with the process $\lambda^{(C)}$ instead of $\lambda$.
Define a sequence of associated random variables $F_n(C)$ similar to $F_n$ as well but with the process $\lambda^{(C)}$ replacing $\lambda$, $\{\tau_i^{n,C}:i \geq 0\}$ replacing $\{\tau_i^{n}:i \geq 0\}$ and ${N}_n^C(T)$ replacing $N_n(T)$, and likewise for $F(C)$ on $(\widetilde{\Omega},\widetilde{\mathcal{G}},\widetilde{\mathbb{P}})$. If for each $C>1$ it holds that
	\begin{align}\label{ass:conv}
		F_n(C)\tols F(C)
	\end{align}  
	and if furthermore
	\begin{align}\label{loc:prelim}
	F_n(C)\1_{\{E_C>T \}} = F_n \1_{\{E_C>T \}} \text{~~ and ~~} F(C)\1_{\{E_C>T \}} = F \1_{\{E_C>T \}}
	\end{align} 
	then we have $F_n\tols F$. 
\end{lemma}

\begin{proof}
Let $\widetilde{\E}$ be the expectation w.r.t.\ $\widetilde{\mathbb{P}}$. We need to prove 
\begin{align*}
&\limsup_{n \to \infty}\left| \E\left[Y f(F_n)\right] - \widetilde{\E}\left[Yf(F)\right] \right| = 0
\end{align*}
where $Y$ is any bounded random variable on $(\Omega,\mathcal{G})$ and $f$ is any bounded continuous function, and using
\begin{align*}
&\limsup_{n \to \infty}\left| \E\left[Y f(F_n)\right] - \widetilde{\E}\left[Yf(F)\right] \right| = \limsup_{C \to \infty} \limsup_{n \to \infty} \left| \E\left[Y f(F_n)\right] - \widetilde{\E}\left[Yf(F)\right] \right| \\ \le& \limsup_{C \to \infty} \limsup_{n \to \infty}\left| \E\left[Y f(F_n)\right] - \E\left[Yf(F_n(C))\right] \right| \\ +&  \limsup_{C \to \infty} \limsup_{n \to \infty}\left| \E\left[Y f(F_n(C))\right] - \widetilde{\E}\left[Yf(F(C))\right] \right| \\ +& \limsup_{C \to \infty} \limsup_{n \to \infty} \left| \widetilde{\E}\left[Y f(F(C))\right] - \widetilde{\E}\left[Yf(F)\right] \right| 
\end{align*} 
it is sufficient to prove that each of the three summands vanishes. For the first one, by boundedness of $Y$ and $f$ and using \eqref{loc:prelim}, it is obvious that 
\[
\left|\E\left[Y (f(F_n)-f(F_n(C)))\right] \right| = \left|\E\left[Y (f(F_n)-f(F_n(C))) \1_{\{ E_C\leq T \}}\right] \right| \le K \mathbb{P}\left(E_C\leq T\right).
\]
Here and below, $K$ always denotes a generic positive constant. Thus,
\[
\limsup_{C \to \infty} \limsup_{n \to \infty}\left| \E\left[Y f(F_n)\right] - \E\left[Yf(F_n(C))\right] \right| \le K \limsup_{C \to \infty} \mathbb{P}\left(E_C\leq T\right) = 0,
\] 
and the same proof applies for the third term. Finally, note that
\[
\limsup_{n \to \infty}\left| \E\left[Y f(F_n(C))\right] - \widetilde{\E}\left[Yf(F(C))\right] \right| = 0
\]
for each fixed $C$ is an immediate consequence of (\ref{ass:conv}). 
\end{proof}

By construction $\lambda_t$ and $\lambda^{(C)}_t$ coincide on the set $\{E_C\leq T\}$ for all $0\leq t \leq T$. As our estimators only deal with observations up to a fixed time horizon $T$ (in our specific case the convenient but arbitrary $T=1$) it is clear that condition \eqref{loc:prelim} is indeed met. Therefore we may assume  for the following proofs that part (a) of Condition \ref{ass:strongstopping} is in force and only prove \eqref{ass:conv} under this strengthened assumption.

Finally, we need to explain why we can assume that part b) of Condition \ref{ass:strongstopping} holds as well. Here we refer to part 2) of the proof of Lemma 9 in \cite{Jacod2018} where a family of discretisation schemes with the desired properties is constructed and where each member of the family coincides with the original sampling scheme up to some random time $S^n_{\ell_n}$. As it is shown that these times converge to infinity almost surely, the same argument as before allows to assume part b) of Condition \ref{ass:strongstopping} without loss of generality. 

For further information on random discretisation schemes one can consult Section 14.1 in \cite{discret} where a slightly different version of Lemma \ref{Lem:Stop1} and other important properties of objects connected to these schemes are proven. We want to name one of those properties in particular because we will use it repeatedly in the following chapters: (14.1.10) in \cite{discret} proves that for all $t\geq 0$ we have
\begin{align}\label{conv:Nn}
\De_n \Nn \pn \int_{0}^{1}\frac{1}{\lambda_s}ds
\end{align}
which basically allows us to treat the random $\Nn$ like the deterministic $\De_n^{-1}$ in all asymptotic considerations.

\subsection{A crucial decomposition}
The proofs of Theorems \ref{thmlncon} and \ref{thmlnclt} rely on a simple decomposition which allows us to identify the terms that play a dominant role in the asymptotic treatment. Precisely, we have 
\begin{align} \label{eq:decomp}
\widetilde{L}^n(p,u)-L(p,u,\beta)=\frac{1}{\Nn-k_n-2} \left(R_1^n(u)+R_2^n(u)+Z^n(u)+R_3^n(u)+R_4^n(u)\right)
\end{align}	
where
\begin{align*}
&z_i(u)=\cos\left(u\frac{\sigma_{\tau_{i-2}^n}\left(\widetilde{\Delta_{i}^nS}-\left(\frac{\lambda_{\tau_{i-2}^n}}{\lambda_{\tau_{i-3}^n}}\right)^{\frac{1}{\beta}-1}\widetilde{\Delta_{i-1}^nS}\right)}{\widetilde{V}_i^n(p)^{1/p}}\right)\\
&~~~~~~~~~~-\E_{i-2}^n\left[\exp\left(-\frac{A_\beta u^\beta{|\sigma_{\tau_{i-2}^n}|}^\beta|\lambda_{\tau_{i-2}^n}|^{1-\beta} ((\phi_{i}^n)^{1-\beta}+(\phi_{i-1}^n)^{1-\beta})}{\Delta_n^{-1}\widetilde{V}_i^n(p)^{\beta/p}}\right)\right],\\
&Z^n(u)=\sum_{i=k_n+3}^{\Nn}z_i(u),
\end{align*}
drives the asymptotics while the residual terms are given by
\begin{align*}
&r_i^1(u)=\cos\left(u\frac{\widetilde{\Delta_i^nX}-\widetilde{\Delta_{i-1}^nX}}{\widetilde{V}_i^n(p)^{1/p}}\right)-\cos\left(u\frac{\sigma_{\tau_{i-2}^n}(\widetilde{\Delta_{i}^nS}-\widetilde{\Delta_{i-1}^nS})}{\widetilde{V}_i^n(p)^{1/p}}\right),\\
&r_i^2(u)=\cos\left(u\frac{\sigma_{\tau_{i-2}^n}(\widetilde{\Delta_{i}^nS}-\widetilde{\Delta_{i-1}^nS})}{\widetilde{V}_i^n(p)^{1/p}}\right)-\cos\left(u\frac{\sigma_{\tau_{i-2}^n}\left(\widetilde{\Delta_{i}^nS}-\left(\frac{\lambda_{\tau_{i-2}^n}}{\lambda_{\tau_{i-3}^n}}\right)^{\frac{1}{\beta}-1}\widetilde{\Delta_{i-1}^nS}\right)}{\widetilde{V}_i^n(p)^{1/p}}\right)\\
&r_i^3(u)=\E_{i-2}^n\left[\exp\left(-\frac{A_\beta u^\beta{|\sigma_{\tau_{i-2}^n}|}^\beta|\lambda_{\tau_{i-2}^n}|^{1-\beta} ((\phi_{i}^n)^{1-\beta}+(\phi_{i-1}^n)^{1-\beta})}{\Delta_n^{-1}\widetilde{V}_i^n(p)^{\beta/p}}\right)\right]\\
&~~~~~~~~~~-\E_{i-2}^n\left[\exp\left(-\frac{C_{p,\beta} u^\beta{|\sigma_{\tau_{i-2}^n}|}^\beta|\lambda_{\tau_{i-2}^n}|^{1-\beta} ((\phi_{i}^n)^{1-\beta}+(\phi_{i-1}^n)^{1-\beta})}{(|\overline{\sigma\lambda}|_i^p)^{\beta/p}}\right)\right],\\
&r_i^4(u)=\E_{i-2}^n\left[\exp\left(-\frac{C_{p,\beta} u^\beta{|\sigma_{\tau_{i-2}^n}|}^\beta|\lambda_{\tau_{i-2}^n}|^{1-\beta} ((\phi_{i}^n)^{1-\beta}+(\phi_{i-1}^n)^{1-\beta})}{(|\overline{\sigma\lambda}|_i^p)^{\beta/p}}\right)\right]-L(p,u,\beta),\\
&R_j^n(u)=\sum_{i=k_n+3}^{\Nn}r_i^j(u) \text{ for } j\in\{1,2,3,4\}.
\end{align*}
Here we have set
\begin{align*}
|\overline{\sigma\lambda}|_i^p=\frac{1}{k_n}\sum_{j=i-k_n-1}^{i-2}|\sigma_{\tau_{j-2}^n}|^p|\lambda_{\tau_{j-2}^n}|^{\frac{p}{\beta}-p}
\end{align*}
and we use the short hand notation $\E_i^n[\cdot]$ in place of $\E[\cdot|\fF^n_{\tau_i^n}]$. We also introduce the notation
\begin{align*}
&\overline{V}_i^n(p)=\frac{1}{k_n}\sum_{j=i-k_n-1}^{i-2}\mathbb{E}^n_{j-2}\left[|\widetilde{\Delta_j^nX}-\widetilde{\Delta_{j-1}^nX}|^p\right], \quad k_n+3\leq i \leq \Nn,
\end{align*}
and set 
\begin{align*}
&\overline{z}_i(u)=\cos\left(u\frac{\lambda_{\tau_{i-2}^n}^{1-1/\beta}\widetilde{\Delta_i^n S}-\lambda_{\tau_{i-3}^n}^{1-1/\beta}\widetilde{\Delta_{i-1}^nS}}{\Delta_n^{1/\beta}\mu_{p,\beta}^{1/\beta}\kappa_{p,\beta}^{1/\beta}}\right)-L(p,u,\beta), \qquad \overline{Z}^n(u)=\sum_{i=k_n+3}^{\Nn}\overline{z}_i(u).
\end{align*}

We will start with a discussion of the asymptotic orders of the residuals for which we always assume that Conditions \ref{ass:strong} and \ref{ass:strongstopping} as well as $k_n\sim C_1 \De_n^{-\varpi}$ for some $C_1> 0$ and some $\varpi\in(0,1)$ are in place. Naturally, we need some preparation to obtain asymptotic negligibility and we will start with a lemma containing a series of bounds for moments of certain increments of (often integrated and rescaled) processes. We will not give a proof of this result but refer to \cite{todorov2015} and \cite{todorov2017}. In fact, the techniques used for the proof will for most parts resemble the ones given therein. The main difference is that our arguments often involve the additional process $\la$ which sometimes complicates matters considerably. 

Recall the decomposition \eqref{decompL}, and we further write 
\begin{align} \label{decompS}
S_t=S_t^{(1)}+S_t^{(2)}
\end{align} 
with $S_t^{(1)}=\int_{0}^{t}\int_\R \kappa(x)\widetilde{\mu}(ds,dx)$, where $\widetilde{\mu}(ds,dx)$ denotes the compensated jump measure of $S$, and $S_t^{(2)}=\int_{0}^{t}\int_\R \kappa'(x){\mu}(ds,dx)$. 

\begin{lemma} \label{lemaux}
Let $i \ge 2$ be arbitrary. 
\begin{itemize}
	\item[(a)] For every $p \in (-1,\be)$ we have
	\begin{align*}
	&\E_{i-2}^n\left[\left|\Delta_n^{-1/\beta}(\widetilde{\Delta_{i}^nS}-\widetilde{\Delta_{i-1}^nS)}\right|^p\right] \le K, \\
	&\E_{i-2}^n\left[\left|\Delta_n^{-1/\beta}(\lambda_{\tau_{i-2}^n}^{1-1/\beta}\widetilde{\Delta_{i}^nS}-\lambda_{\tau_{i-3}^n}^{1-1/\beta}\widetilde{\Delta_{i-1}^nS)}\right|^p\right] = \kappa_{p,\beta}^{p/\beta}\mu_{p,\beta}^{p/\beta},
	\end{align*}
	with the constants from (\ref{defcon}). 
	\item[(b)] For every $p\in(0,\beta)$ we have
	\begin{align*}
	\left|\E_{i-2}^n\left[\left|\Delta_n^{-1/\beta}(\widetilde{\Delta_{i}^nS}-\widetilde{\Delta_{i-1}^nS)}\right|^p\right] - \lambda_{\tau_{i-2}^n}^{p/\beta - p}\kappa_{p,\beta}^{p/\beta}\mu_{p,\beta}^{p/\beta} \right| \leq K \Delta_n^{1/2}.
	\end{align*}
	\item[(c)] For every $q\in (0,2]$, every $\iota > 0$ and every $l \in \{0,1\}$ we have
	\begin{align*}
	\E_{i-2}^n\left[\left|\frac{\Delta_n}{\tau_{i-2+l}^n-\tau_{i-1+l}^n}\int_{\tau_{i-2+l}^n}^{\tau_{i-1+l}^n}(\sigma_{u-}-\sigma_{\tau_{i-2}^n})dS_u^{(1)}\right|^q\right]	 \leq K \Delta_n^{q/2+q/\beta\wedge 1-\iota}.
	\end{align*}
		\item[(d)] For every $q\in (0,2]$ and every $l \in \{0,1\}$ we have
	\begin{align*}
	\E_{i-2}^n\left[\left|\frac{\Delta_n}{\tau_{i-1+l}^n-\tau_{i-2+l}^n}\int_{\tau_{i-2+l}^n}^{\tau_{i-1+l}^n}(\sigma_{u-}-\sigma_{\tau_{i-2}^n})dS_u^{(2)}\right|^q\right] \leq K \Delta_n^{q/2+1}.
	\end{align*}
	\item[(e)] For every $q> 0$, every $\iota > 0$ and every $l \in \{0,1\}$ we have
	\begin{align*}
	\E_{i-2}^n\left[\left|\frac{\Delta_n}{\tau_{i-1+l}^n-\tau_{i-2+l}^n}\int_{\tau_{i-2+l}^n}^{\tau_{i-1+l}^n}(\sigma_{u-} -\sigma_{\tau_{i-2}^n}) d\acute{S}_u\right|^q\right] \leq K\Delta_n^{q/2+{(q/\beta')\wedge 1}-\iota},
	\end{align*}
	and the same relation holds with $\grave{S}$ instead of $\acute{S}$.
	\item[(f)] For every $q > 0$ and every $l \in \{0,1\}$ we have
	\begin{align*}
	\E_{i-2}^n\left[\left|\frac{\Delta_n}{\tau_{i-1+l}^n-\tau_{i-2+l}^n}\Delta_{i-1+l}^nY\right|^q\right]\leq K\Delta_n^{(q/\beta')\wedge 1}.
	\end{align*}
	\item[(g)] For every $q\in (0,2]$ we have
	\begin{align*}
	\E_{i-2}^n\left[\left|\frac{\Delta_n}{\tau_{i}^n-\tau_{i-1}^n}\int_{\tau_{i-1}^n}^{\tau_{i}^n}\alpha_{u}du-\frac{\Delta_n}{\tau_{i-1}^n-\tau_{i-2}^n}\int_{\tau_{i-2}^n}^{\tau_{i-1}^n}\alpha_{u}du\right|^q\right] \leq \Delta_n^{3q/2}.
	\end{align*}
\end{itemize}
\end{lemma}

A second lemma, again without proof, discusses bounds for moments of increments of semimartingales. Again, it has some resemblance to results in  \cite{todorov2015} and \cite{todorov2017} but its proof is slightly more involved due to the random observation scheme.

\begin{lemma}\label{lem:ineq_power}
	Let $A$ be a semimartingale satisfying the same properties as $\si$ in Assumption \ref{ass:strong}. Then, for any $-1<p<1$ and any $y>0$ and with $K$ possibly depending on $y$ we have 
	\begin{align*}
	\left|\E_i^n\left[|A_{\tau^n_{i+j}}|^p-|A_{\tau^n_{i}}|^p\right]\right| &\leq K j \Delta_n, \\
	\E_i^n\left[\left||A_{\tau^n_{i+j}}|^p-|A_{\tau^n_{i}}|^p\right|^y\right] &\leq K (j\Delta_n)^{y/2\wedge1}.
	\end{align*}
\end{lemma}

After the presentation of these auxiliary claims, we focus on results which directly simplify the discussion of the asymptotic negligibility of the residual terms. We begin with a results that helps in the treatment of $R_1^n$.  

\begin{lemma} \label{lemaln}
Let $\iota>0$ and $0<p<\frac{\beta}{2}$ be arbitrary. Then, for any $i \ge 2$, we have 
	\begin{align*} 
	\Delta_n^{-p/\beta}\mathbb{E}_{i-2}^n\left[\left||\widetilde{\Delta_i^nX} -\widetilde{\Delta_{i-1}^nX}|^p-|\sigma_{\tau_{i-2}^n}|^p|\widetilde{\Delta_{i}^nS}-\widetilde{\Delta_{i-1}^nS}|^p\right|\right]\leq K\alpha_n
	\end{align*}
	with $\alpha_n=\Delta_n^{\frac{\beta}{2}\frac{p+1}{\beta+1}\wedge((\frac{p}{\beta'}\wedge 1)-\frac{p}{\beta})\wedge\frac{1}{2}-\iota}$.
\end{lemma} 

\begin{proof}
The proof relies on bounds for moments of several stochastic integrals, mostly connected with the jump process $L_t$ and its parts. Using (\ref{decompL}) and (\ref{decompS}) we may write $\widetilde{\Delta_i^nX} -\widetilde{\Delta_{i-1}^nX} = \chi_i^{(n,1)} + \chi_i^{(n,2)} + \chi_i^{(n,3)}$ with 
\begin{align*}
	&\chi_i^{(n,1)}=\sigma_{\tau_{i-2}^n}\left( \widetilde{\Delta_i^n S} -\widetilde{\Delta_{i-1}^n S} \right),\\
	&\chi_i^{(n,2)} =\frac{\Delta_n}{\tau_{i}^n-\tau_{i-1}^n}\int_{\tau_{i-1}^n}^{\tau_{i}^n}(\sigma_{u-}-\sigma_{\tau_{i-2}^n})dS_u^{(1)}-\frac{\Delta_n}{\tau_{i-1}^n-\tau_{i-2}^n}\int_{\tau_{i-2}^n}^{\tau_{i-1}^n}(\sigma_{u-}-\sigma_{\tau_{i-2}^n})dS_u^{(1)}\\	
	&~~~~~~~~+\frac{\Delta_n}{\tau_{i}^n-\tau_{i-1}^n}\int_{\tau_{i-1}^n}^{\tau_{i}^n}\alpha_{u}du-\frac{\Delta_n}{\tau_{i-1}^n-\tau_{i-2}^n}\int_{\tau_{i-2}^n}^{\tau_{i-1}^n}\alpha_{u}du,\\
	& \chi_i^{(n,3)}=\frac{\Delta_n}{\tau_{i}^n-\tau_{i-1}^n}\int_{\tau_{i-1}^n}^{\tau_{i}^n}(\sigma_{u-}-\sigma_{\tau_{i-2}^n})dS_u^{(2)}-\frac{\Delta_n}{\tau_{i-1}^n-\tau_{i-2}^n}\int_{\tau_{i-2}^n}^{\tau_{i-1}^n}(\sigma_{u-}-\sigma_{\tau_{i-2}^n})dS_u^{(2)}\\
	&~~~~~~~~+\frac{\Delta_n}{\tau_{i}^n-\tau_{i-1}^n}\Delta_i^nY-\frac{\Delta_n}{\tau_{i-1}^n-\tau_{i-2}^n}\Delta_{i-1}^nY\\ &~~~~~~~~+\frac{\Delta_n}{\tau_{i}^n-\tau_{i-1}^n}\int_{\tau_{i-1}^n}^{\tau_{i}^n}\sigma_{u-}d\acute{S}_u -\frac{\Delta_n}{\tau_{i-1}^n-\tau_{i-2}^n}\int_{\tau_{i-2}^n}^{\tau_{i-1}^n}\sigma_{u-}d\acute{S}_u \\
	&~~~~~~~~-\frac{\Delta_n}{\tau_{i}^n-\tau_{i-1}^n}\int_{\tau_{i-1}^n}^{\tau_{i}^n}\sigma_{u-}d\grave{S}_u +\frac{\Delta_n}{\tau_{i-1}^n-\tau_{i-2}^n}\int_{\tau_{i-2}^n}^{\tau_{i-1}^n}\sigma_{u-}d\grave{S}_u.
	\end{align*}
We obviously have 
\begin{align*}
&\Delta_n^{-p/\beta} \mathbb{E}_{i-2}^n\left[\left||\chi_i^{(n,1)} + \chi_i^{(n,2)} + \chi_i^{(n,3)}|^p - |\chi_i^{(n,1)}|^p \right| \right] \\ \leq& \Delta_n^{-p/\beta} \mathbb{E}_{i-2}^n\left[\left||\chi_i^{(n,1)} + \chi_i^{(n,2)} + \chi_i^{(n,3)}|^p - |\chi_i^{(n,1)} + \chi_i^{(n,2)}|^p \right| \right] \\ +& \Delta_n^{-p/\beta} \mathbb{E}_{i-2}^n\left[\left||\chi_i^{(n,1)} + \chi_i^{(n,2)}|^p - |\chi_i^{(n,1)}|^p \right| \right],
\end{align*}
and since $p < \be/2 < 1$ holds, the inequality 
\begin{align} \label{ineqchi3}
\Delta_n^{-p/\beta} \mathbb{E}_{i-2}^n\left[|\chi_i^{(n,3)}|^p\right] \le K \Delta_n^{p/\beta'\wedge 1-p/\beta},
\end{align}
which is any easy consequence of parts (d), (e) and (f) of Lemma \ref{lemaux}, is enough to fully treat the first term on the right hand side.
For the second term we will use parts (a), (c) and (g) of Lemma \ref{lemaux} plus the algebraic inequality 
	\begin{align*}
	\left||a+b|^p-|a|^p\right|\leq K|a|^{p-1}|b|\1_{\{|a|>\epsilon,|b|<\frac{1}{2}\epsilon\}}+|b|^p(\1_{\{|a|\leq \epsilon\}}+\1_{\{|b|> \frac{1}{2}\epsilon\}}),
	\end{align*}
	which holds for any $\epsilon>0$ and $p\in(0,1]$ and a constant $K$ that does not depend on $\epsilon$ and which we apply with $a = \Delta_n^{-1/\beta} \chi_i^{(n,1)}$ and $b=\Delta_n^{-1/\beta} \chi_i^{(n,2)}$. We start with the latter two terms and let $0 < \epsilon < 1$ be arbitrary. Then Markov inequality in combination with Hölder inequality first gives 
	\begin{align*}
	&\E_{i-2}^n\left[|\Delta_n^{-1/\beta}\chi_i^{(n,2)}|^p\1_{\{|\Delta_n^{-1/\beta} \chi_i^{(n,1)}|\leq \epsilon\}}\right]\\ \leq& \left(\E_{i-2}^n\left[|\Delta_n^{-1/\beta} \chi_i^{(n,1)}|^{-1+\iota}\right]\epsilon^{1-\iota}\right)^{1-\frac{p}{\beta}} \left(\E_{i-2}^n\left[|\Delta_n^{-1/\beta} \chi_i^{(n,2)}|^\beta\right]\right)^\frac{p}{\beta} \\ \leq& K \epsilon^{1-p/\beta-\iota} (\Delta_n^{\beta/2-\iota})^{p/\beta}	\leq K \epsilon^{1-p/\beta-\iota}\Delta_n^{p/2-\iota}
	\end{align*}
(with a slight abuse of notation but remember that $\iota > 0$ can be chosen arbitrarily)	and then
	\begin{align*}
	&\E_{i-2}^n\left[|\Delta_n^{-1/\beta}\chi_i^{(n,2)}|^p\1_{\{|\Delta_n^{-1/\beta}\chi_i^{(n,2)}|> \frac{1}{2}\epsilon\}}\right] \\ \leq& \left(\frac{\E_{i-2}^n\left[|\Delta_n^{-1/\beta}\chi_i^{(n,2)}|^\beta\right]}{(\frac{1}{2}\epsilon)^\beta}\right)^{1-\frac{p}{\beta}} \left(\E_{i-2}^n\left[|\Delta_n^{-1/\beta}\chi_i^{(n,2)}|^\beta\right]\right)^{\frac{p}{\beta}} \leq K \epsilon^{-(\beta-p)} \Delta_n^{\beta/2-\iota}
	\end{align*}
	as well. Setting $\epsilon=\Delta_n^{\frac{1}{2}\frac{\beta}{\beta+1}}$ gives $K \Delta_n^{\frac{\beta}{2}\frac{p+1}{\beta+1}-\iota}$ as the upper bound in both terms.
	
	Finally, we have to distinguish between $p > 1/\be$ and $p \le 1/\be$. In the first case a simple application of Hölder inequality gives
	\begin{align*}
	\Delta_n^{-p/\beta} \mathbb{E}_{i-2}^n\left[|\chi_i^{(n,1)}|^{p-1} |\chi_i^{(n,2)}| \right] &\leq  \Delta_n^{-p/\beta}  \left(\E_{i-2}^n\left[|\chi_i^{(n,1)}|^{\frac{(p-1)\be}{\be-1}}\right]\right)^{\frac{\be-1}{\be}}\left(\E_{i-2}^n| \chi_i^{(n,2)}|^\be\right)^{\frac{1}{\be}} \\ &\leq K \Delta_n^{-1/\be + 1/2+1/\beta-\iota} = K \Delta_n^{1/2-\iota}
	\end{align*}
	for our specific choice of $\iota > 0$. Note that (a) in Lemma \ref{lemaux} was indeed applicable as $p > 1/\be$ ensures $(p-1)\be/(\be-1) > -1$. In the second case we set $\epsilon$ as above and use Markov inequality with $r = \frac{1+\iota}{\be} - p$. Then 
	\begin{align*}
	&\Delta_n^{-p/\beta} \mathbb{E}_{i-2}^n\left[|\chi_i^{(n,1)}|^{p-1} |\chi_i^{(n,2)}| \1_{\{|\Delta_n^{-1/\beta} \chi_i^{(n,1)}|> \epsilon\}}\right] \\ \leq&  \Delta_n^{-(p+r)/\beta} \mathbb{E}_{i-2}^n\left[|\chi_i^{(n,1)}|^{p+r-1} |\chi_i^{(n,2)}| \right] \epsilon^{-r},
	\end{align*}
	and as now $p+r > 1/\beta$ by construction, the same proof as in the first case proves this term to be of the order $\De_n^{1/2-\iota} \epsilon^{-r}$. Then 
	\[
	\frac 12-\iota - r \frac{\beta}{2(\beta+1)} = \frac 12- \frac{1}{2(\beta +1)} + \frac{p\be}{2(\beta +1)} -\iota = \frac{\beta}{2}\frac{p+1}{\beta+1}-\iota
	\]
	(with the same abuse of notation) ends the proof. 
\end{proof}

We also need to control the denominators in $R_1^n$ and $R_2^n$ to make sure that they are bounded away from zero with high probability. To this end, we need two auxiliary results, and in both cases we let $i \ge k_n + 3$ and $0<p<\frac{\beta}{2}$ be arbitrary. The first result deals with the variables $\overline{V}_i^n(p)$ which we introduced before, and it will be used for the treatment of $R_3^n$ later on as well. Its proof is omitted as it is essentially the same as the one for equation (9.4) in \cite{todorov2015} and exploits standard inequalities for discrete martingales.

\begin{lemma} \label{lemoverline}
	Let $k_n \sim C_1 \De_n^{-\varpi}$ for some $C_1 > 0$ and $\varpi \in (0,1)$. Then for all $1\leq x<\frac{\beta}{p}$ we have
	\begin{align*}
	\Delta_n^{-xp/\beta}\E \left[{|\widetilde{V}_i^n(p)-\overline{V}_i^n(p)|}^x \right]\leq K_xk_n^{-x/2}
	\end{align*}
	where the constant $K_x$ might depend on $x$. 
\end{lemma}	

The second result deals with the set $$\mathcal{C}_i^n=\left \{|\Delta_n^{-p/\beta}\widetilde{V}_i^n(p)-\sigs_i\mu_{p,\beta}^{p/\beta}\kappa_{p,\beta}^{p/\beta}|>\frac{1}{2}\sigs_i\mu_{p,\beta}^{p/\beta}\kappa_{p,\beta}^{p/\beta} \right \}$$ and bounds its probability. 

\begin{lemma} \label{lemC}
	For every fixed $\iota > 0$ we have
	\begin{align*}
	\mathbb{P}(\mathcal{C}_i^n)\leq K_\iota k_n^{-\beta/2p+\iota}
	\end{align*}
	where the constant $K_\iota$ might depend on $\iota$. 
\end{lemma}

\begin{proof}
Lemmas \ref{lemaux} and \ref{lemaln} give
\begin{align}
	\nonumber &\left|\Delta_n^{-p/\beta}\E_{i-2}^n\left[\left|\widetilde{\Delta_i^nX}-\widetilde{\Delta_{i-1}^nX}\right|^p\right]-|\sigma_{\tau_{i-2}^n}|^p|\lambda_{\tau_{j-2}^n}|^{\frac{p}{\beta}-p}\mu_{p,\beta}^{p/\beta}\kappa_{p,\beta}^{p/\beta}\right|\\
	\nonumber &\leq\Delta_n^{-p/\beta}\E_{i-2}^n\left[\left|\left|\widetilde{\Delta_i^nX}-\widetilde{\Delta_{i-1}^nX}\right|^p- |\sigma_{\tau_{i-2}^n}|^p\left|\widetilde{\Delta_i^nS}-\widetilde{\Delta_{i-1}^nS}\right|^p\right|\right]\\
	\nonumber&~~~+\left|\E_{i-2}\left[\Delta_n^{-p/\beta}|\sigma_{\tau_{i-2}^n}|^p\left|\widetilde{\Delta_i^nS}-\widetilde{\Delta_{i-1}^nS}\right|^p\right]-|\sigma_{\tau_{i-2}^n}|^p|\lambda_{\tau_{j-2}^n}|^{\frac{p}{\beta}-p}\mu_{p,\beta}^{p/\beta}\kappa_{p,\beta}^{p/\beta}\right|\\
	&\leq K\alpha_n + K |\sigma_{\tau_{i-2}^n}|^p \Delta_n^{1/2},  \label{ineqas} 
	\end{align}
so 
\begin{align} \label{ineqoverline}
	|\Delta_n^{-p/\beta}\overline{V}_i^n(p)-\sigs_i\mup^{p/\beta}\kappa_{p,\beta}^{p/\beta}| \leq K\alpha_n + K \sigl_i \Delta_n^{1/2},
	\end{align}
with $\sigl_i$ being defined as $\sigs_i$ but with $\la = 1$. Now, it is a simple consequence of Conditions \ref{ass:strong} and \ref{ass:strongstopping} that both $\sigs_i$ and $\sigl_i$ are uniformly bounded from above and below. Using $\alpha_n\rightarrow 0$ and $\Delta_n^{1/2}\rightarrow 0$ there then exists some $n_0\in \mathbb{N}$ such that 
	\begin{align*}
	\mathbb{P}\left( \left|\Delta_n^{-p/\beta}\overline{V}_i^n(p)-\sigs_i\mu_{p,\beta}^{p/\beta}\kappa_{p,\beta}^{p/\beta}\right|>\frac{1}{4}\sigs_i\mu_{p,\beta}^{p/\beta}\kappa_{p,\beta}^{p/\beta} \right) = 0
	\end{align*}
 for all $n\geq n_0$. For these $n$,
\begin{align*}
	\mathbb{P}(\Ci)&\leq \mathbb{P}\left(\Delta_n^{-p/\beta}|\widetilde{V}_i^n(p)-\overline{V}_i^n(p)|>\frac{1}{4}\sigs_i\mu_{p,\beta}^{p/\beta}\kappa_{p,\beta}^{p/\beta}\right)\\
	&~~~~+\mathbb{P}\left(|\Delta_n^{-p/\beta}\overline{V}_i^n(p)-\sigs_i\mu_{p,\beta}^{p/\beta}\kappa_{p,\beta}^{p/\beta}|>\frac{1}{4}\sigs_i\mu_{p,\beta}^{p/\beta}\kappa_{p,\beta}^{p/\beta}\right) \leq K_\iota k_n^{-\beta/2p+\iota}
	\end{align*}
from Lemma \ref{lemoverline}, with a choice of $x$ arbitrarily close to $\frac{\beta}{p}$. The claim follows.
\end{proof}

Finally, we provide two lemmas that simplify the discussion of $R_4^n$. The first one gives an alternative representation for the limiting variable $L(p,u,\beta)$. 

\begin{lemma}\label{lem:exp}
	It holds that
	\begin{align*}
	\E^n_{i-2}\left[\cos\left(u\frac{\lambda_{\tau_{i-2}^n}^{1-1/\beta}\widetilde{\Delta_i^n S}-\lambda_{\tau_{i-3}^n}^{1-1/\beta}\widetilde{\Delta_{i-1}^nS}}{\Delta_n^{1/\beta}\mu_{p,\beta}^{1/\beta}\kappa_{p,\beta}^{1/\beta}}\right)\right]=L(p,u,\beta).
	\end{align*}
\end{lemma} 

\begin{proof}
For the sake of simplicity, let us assume that the probability space can be even further enlarged to accomodate three independent random variables $S^{(1)}$, $S^{(2)}$ and $S^{(3)}$, independent of $\gG$, all with the same distribution as $S_1$, i.e.\ distributed as a L\'evy process with characteristic triplet $(0,0,F)$ at time 1, $F(dx) = F(x) dx$ with $F(x) =  {A}|x|^{-(1+\beta)}$. Using standard properties of stable processes (see e.g.\ Section 1.2 in \cite{stableprocesses1994}), for constants $\sigma_1, \sigma_2 \in \R$ we have
	\begin{align*}
	\sigma_1 S^{(1)} + \sigma_2 S^{(2)} \sim (|\sigma_1|^\beta+|\sigma_2|^\beta)^{1/\beta} S^{(3)},
	\end{align*}
	and for our original process $S_t$ the stability relation 
	\begin{align*}
	&(S_{t+r}-S_{t})\sim r^{1/\beta}S_1 \text{ for all } r, t \ge 0
	\end{align*}
	holds as well. 
	Because the increments of the process $(S_t)_{t\geq \tau_{i-2}^n}$ are independent of $\tau_{i}^n-\tau_{i-1}^n= \Delta_n \phi_i^n \lambda_{\tau_{i-2}^n}$ the conditional distribution of $\Delta_{i}^nS$ given $(\tau_{i}^n-\tau_{i-1}^n)=a$ equals the one of $a^{1/\beta}S^{(1)}$.	Thus, for all Borel sets $M$ we obtain e.g.\
	\begin{align*}
	\E\left[\1_M\left(\frac{\Delta_{i}^nS}{\tau_{i}^n-\tau_{i-1}^n}\right)\right]&=\int_\R \E\left[\1_M\left(\frac{\Delta_{i}^nS}{\tau_{i}^n-\tau_{i-1}^n}\right)\Bigg\vert\left(\tau_{i}^n-\tau_{i-1}^n\right)=a\right]\mathbb{P}^{(\tau_{i}^n-\tau_{i-1}^n)}(da)\\
	&=\E\left[\1_M\left((\tau_{i}^n-\tau_{i-1}^n)^{1/\beta-1}\right)S^{(1)}\right]
	\end{align*}
	using that we assume that all moments of $(\phi_i^n)^{q}$ for $q\in(-2,0)$ exist, as well as $1/\beta-1>-1$. Put differently,
	\begin{align*}
	&\Delta_n^{-1/\beta}\widetilde{\Delta_{i}^nS}\sim \Delta_n^{1-1/\beta} (\tau_{i}^n-\tau_{i-1}^n)^{1/\beta-1}S^{(1)}\sim (\lambda_{\tau_{i-2}^n}\phi_i^n)^{1/\beta-1}S^{(1)}, \\
	&\Delta_n^{-1/\beta}\widetilde{\Delta_{i-1}^nS}\sim \Delta_n^{1-1/\beta} (\tau_{i-1}^n-\tau_{i-2}^n)^{1/\beta-1}S^{(2)}\sim (\lambda_{\tau_{i-3}^n}\phi_{i-1}^n)^{1/\beta-1}S^{(2)},
	\end{align*}
	where $\phi_{i-1}^n$, $\phi_i^n$, $S^{(1)}$, $S^{(2)}$ are all independent of $\mathcal{F}_{\tau_{i-2}^n}$ and of each other. Thus
		\begin{align}
	&\E^n_{i-2}\left[\cos\left(u\frac{\lambda_{\tau_{i-2}^n}^{1-1/\beta}\widetilde{\Delta_i^n S}-\lambda_{\tau_{i-3}^n}^{1-1/\beta}\widetilde{\Delta_{i-1}^nS}}{\Delta_n^{1/\beta}\mu_{p,\beta}^{1/\beta}\kappa_{p,\beta}^{1/\beta}}\right)\right] \nonumber \\
	=&\E_{i-2}^n\left[\cos\left(u\frac{((\phi_i^n)^{1-\beta}+(\phi_{i-1}^n)^{1-\beta})^{1/\beta}S^{(3)}}{\mu_{p,\beta}^{1/\beta}\kappa_{p,\beta}^{1/\beta}}\right)\right]\nonumber \\
	=&\E\left[\exp\left(-u^\beta\frac{A_{\beta}}{\mu_{p,\beta}\kappa_{p,\beta}}((\phi_{i}^n)^{1-\beta}+(\phi_{i-1}^n)^{1-\beta})\right)\right]\nonumber \\
	=&\E\left[\exp\left(-u^\beta C_{p,\beta}((\phi_{i}^n)^{1-\beta}+(\phi_{i-1}^n)^{1-\beta})\right)\right]=L(p,u,\beta) \label{identL}
	\end{align}
	after successive conditioning.
\end{proof}

Using Lemma \ref{lem:exp} it is clear that the treatment of $R_4^n$ hinges on the question how well $\sigs_i$ can be approximated by $|\sigma_{\tau_{i-2}^n}|^p |\lambda_{\tau_{i-2}^n}|$. 

\begin{lemma} \label{lem:ineqsila}
Let $k_n \sim C_1 \De_n^{-\varpi}$ for some $C_1 > 0$ and $\varpi \in (0,1)$. Then for all $y>1$ we have
	\begin{align*}
	\left|\E_{i-k_n-3}^n\left[\sigs_i-|\sigma_{\tau_{i-2}^n}|^p|\lambda_{\tau_{i-2}^n}|^{\frac{p}{\beta}-p}\right]\right|&\leq K k_n\Delta_n,\\
	\E_{i-k_n-3}\left[\left|\sigs_i-|\sigma_{\tau_{i-2}^n}|^p|\lambda_{\tau_{i-2}^n}|^{\frac{p}{\beta}-p}\right|^y\right]&\leq K(k_n\Delta_n)^{y/2\wedge 1}. 
	\end{align*}
\end{lemma}

\begin{proof}
We start first with a proof of 
	\begin{align*}
	&\left|\E_{i}^n\left[|\sigma_{\tau^n_{i+j}}|^p|\lambda_{\tau_{i+j}^n}|^{\frac{p}{\beta}-p}-|\sigma_{\tau^n_{i}}|^p|\lambda_{\tau_{i}^n}|^{\frac{p}{\beta}-p}\right]\right|\leq K j\Delta_n
	\end{align*}
and
	\begin{align} \label{ineqsila}
	&\E_{i}^n\left[\left||\sigma_{\tau^n_{i+j}}|^p|\lambda_{\tau_{i+j}^n}|^{\frac{p}{\beta}-p}-|\sigma_{\tau^n_{i}}|^p|\lambda_{\tau_{i}^n}|^{\frac{p}{\beta}-p}\right|^y\right]\leq K (j\Delta_n)^{\frac{y}{2}\wedge 1}.
	\end{align}
For the first claim, note that the decomposition 
\begin{align*}
&|\sigma_{\tau^n_{i+j}}|^p|\lambda_{\tau_{i+j}^n}|^{\frac{p}{\beta}-p}-|\sigma_{\tau^n_{i}}|^p|\lambda_{\tau_{i}^n}|^{\frac{p}{\beta}-p} \\=&  |\sigma_{\tau^n_{i}}|^p(|\lambda_{\tau_{i+j}^n}|^{\frac{p}{\beta}-p} - |\lambda_{\tau_{i}^n}|^{\frac{p}{\beta}-p}) + |\lambda_{\tau_{i}^n}|^{\frac{p}{\beta}-p} (|\sigma_{\tau^n_{i+j}}|^p - |\sigma_{\tau^n_{i}}|^p) \\ +& (|\lambda_{\tau_{i+j}^n}|^{\frac{p}{\beta}-p} - |\lambda_{\tau_{i}^n}|^{\frac{p}{\beta}-p}) (|\sigma_{\tau^n_{i+j}}|^p - |\sigma_{\tau^n_{i}}|^p)
\end{align*}
holds. Now, we can apply Lemma \ref{lem:ineq_power} for each of the three terms, for the third one together with Cauchy-Schwarz inequality, and using boundedness of $\si$ and $\la$ from below and above plus the fact that $1 < \be < 2$ and $p < \be/2 < 1$ guarantee the exponents to lie between $-1$ and $1$. A similar reasoning works for the second claim. 

We then obtain 
	\begin{align*}
	&\left|\E_{i-k_n-3}^n[\sigs_i-|\sigma_{\tau_{i-2}^n}|^p|\lambda_{\tau_{i-2}^n}|^{\frac{p}{\beta}-p}]\right|\\
	=&\left|\E_{i-k_n-3}^n\left[\frac{1}{k_n}\sum_{j=i-k_n-1}^{i-2}(|\sigma_{\tau_{j-2}^n}|^p|\lambda_{\tau_{j-2}^n}|^{\frac{p}{\beta}-p}-|\sigma_{\tau_{i-2}^n}|^p|\lambda_{\tau_{i-2}^n}|^{\frac{p}{\beta}-p})\right]\right|\\
	\leq&\frac{1}{k_n}\sum_{j=i-k_n-1}^{i-2}\left|\E_{i-k_n-3}^n\left[|\sigma_{\tau_{j-2}^n}|^p|\lambda_{\tau_{j-2}^n}|^{\frac{p}{\beta}-p}-|\sigma_{\tau_{i-2}^n}|^p|\lambda_{\tau_{i-2}^n}|^{\frac{p}{\beta}-p}\right]\right|\\
	\leq& \frac{1}{k_n}\sum_{j=i-k_n-1}^{i-2}K(i-j)\Delta_n \leq Kk_n\Delta_n
	\end{align*}
	easily, and by convexity of $x\mapsto x^y$ on $\mathbb{R}_+$
	\begin{align*}
	&\E_{i-k_n-3}^n\left[\left|\sigs_i-|\sigma_{\tau_{i-2}^n}|^p|\lambda_{\tau_{i-2}^n}|^{\frac{p}{\beta}-p}\right|^y\right]\\
	=&\E_{i-k_n-3}^n\left[\left|\frac{1}{k_n}\sum_{j=i-k_n-1}^{i-2}(|\sigma_{\tau_{j-2}^n}|^p|\lambda_{\tau_{j-2}^n}|^{\frac{p}{\beta}-p}-|\sigma_{\tau_{i-2}^n}|^p|\lambda_{\tau_{i-2}^n}|^{\frac{p}{\beta}-p})\right|^y\right]\\
	\leq& \frac{1}{k_n}\sum_{j=i-k_n-1}^{i-2}\E_{i-k_n-3}^n\left[\left||\sigma_{\tau_{j-2}^n}|^p|\lambda_{\tau_{i-2}^n}|^{\frac{p}{\beta}-p}-|\sigma_{\tau_{i-2}^n}|^p|\lambda_{\tau_{i-2}^n}|^{\frac{p}{\beta}-p}\right|^y\right]
	\end{align*}
	gives the claim. 
\end{proof}

\subsection{Bounding the residual terms}
In what follows, let $u > 0$, $0<p<\frac{\beta}{2}$ and $\iota > 0$ be arbitrary but fixed, and we always assume that $k_n\sim C_1 \De_n^{-\varpi}$ for some $C_1 > 0$ and $\varpi \in (0,1)$. In the following we also use the notation $n = \De_n^{-1}$ for convenience. 
\begin{lemma}\label{lemma:R1}
	We have
	\begin{align*}
	\frac{1}{n-k_n-2}\mathbb{E}\left[\left|R_1^n(u)\right|\right] \leq K \left(k_n^{-\beta/2p+\iota}\vee  u^{\beta'}\Delta_n^{(\beta-\beta')/\beta} \vee u\Delta_n^{1/2-\iota} \right)
	\end{align*}
	where the constant $K$ might depend on $p$, $\be$ and $\iota$ but not on $u$.
\end{lemma}
\begin{proof}
We decompose $r_i^1(u_n)=r_i^1(u_n) \1_\Ci+r_i^1(u_n)\1_{(\Ci)^C}$, and as $\cos(x)$ is bounded we have for any $i \ge k_n+3$ by Lemma \ref{lemC}
	\begin{align*}
	\E_{i-2}^n\left[\left|r_i^1(u)\1_{\Ci}\right|\right]\leq K \mathbb{P}(\Ci) \leq K k_n^{-\beta/2p+\iota}.
	\end{align*}
Thus, using Assumption \ref{ass:strongstopping} we obtain 
	\begin{align}\label{key:Nn_above}
	&\frac{1}{n-k_n-2}\E\left[\sum_{i=k_n+3}^{\Nn}\left|r_i^1(u) \1_{\Ci}\right|\right] = \frac{1}{n-k_n-2}\E\left[\sum_{i=k_n+3}^{\lfloor C \De_n^{-1} \rfloor}\left|\1_{\{\Nn\geq i\}}r_i^1(u) \1_{\Ci}\right|\right]\nonumber\\
	 \leq& \frac{1}{n-k_n-2}\sum_{i=k_n+3}^{\lfloor C \De_n^{-1} \rfloor}\E\left[|r_i^1(u) \1_{\Ci}|\right] \leq K k_n^{-\beta/2p+\iota}.
	\end{align}
	On the other hand, on $(\Ci)^C$ the relation 
	\[
\frac{1}{2}\sigs_i\mu_{p,\beta}^{p/\beta}\kappa_{p,\beta}^{p/\beta} \le \Delta_n^{-p/\beta}\widetilde{V}_i^n(p)\leq\frac{3}{2}\sigs_i\mu_{p,\beta}^{p/\beta}\kappa_{p,\beta}^{p/\beta}
	\]
	holds. Since $\sigs_i$ is bounded from above and below by Assumption \ref{ass:strong}, $\Delta_n^{-p/\beta}\widetilde{V}_i^n(p)$ is now likewise with a constant possibly depending on $p$ and $\be$. Let us use the notation from the proof of Lemma \ref{lemaln} and write
	\begin{align*}
	\widetilde{\Delta_i^nX} -\widetilde{\Delta_{i-1}^nX} = \chi_i^{(n,1)} + \chi_i^{(n,2)} + \chi_i^{(n,3)}, \quad \sigma_{\tau_{i-2}^n}(\widetilde{\Delta_{i}^nS}-\widetilde{\Delta_{i-1}^nS})=\chi_i^{(n,1)}.
	\end{align*}
	Using the boundedness of $\Delta_n^{-p/\beta}\widetilde{V}_i^n(p)$ on $(\Ci)^C$ and the inequality $|\cos(x)-\cos(y)|\leq2|x-y|^p$ for all $x,y\in\R$ and $p \in(0,1]$ we have
	\begin{align*}
	&\E_{i-2}^n\left[\left|r_i^1(u)\1_{(\Ci)^C}\right|\right]\\ 
	\leq& \E_{i-2}^n\left[\left|\cos\left(u \frac{\chi_i^{(n,1)} + \chi_i^{(n,2)} + \chi_i^{(n,3)}}{(\widetilde{V}_i^n(p))^{1/p}}\right)-\cos\left(u\frac{\chi_i^{(n,1)} + \chi_i^{(n,2)}}{(\widetilde{V}_i^n(p))^{1/p}}\right)\right|\1_{(\Ci)^C}\right]\\+&\E_{i-2}^n\left[\left|\cos\left(u \frac{\chi_i^{(n,1)} + \chi_i^{(n,2)}}{(\widetilde{V}_i^n(p))^{1/p}}\right)-\cos\left(u\frac{\chi_i^{(n,1)}}{(\widetilde{V}_i^n(p))^{1/p}}\right)\right|\1_{(\Ci)^C}\right]\\
	\leq& K \left(\E_{i-2}^n\left[\left|u \Delta_n^{-1/\beta} \chi_i^{(n,3)}\right|^{\beta'}\right]+\E_{i-2}^n\left[\left|u \Delta_n^{-1/\beta} \chi_i^{(n,2)}\right|\right]\right).
	\end{align*}
We then get
	\begin{align*}			
		\E_{i-2}^n\left[\left|u \Delta_n^{-1/\beta} \chi_i^{(n,2)}\right|\right] \leq K u \Delta_n^{1/2-\iota}, \quad
		\E_{i-2}^n\left[\left|u \Delta_n^{-1/\beta} \chi_i^{(n,3)}\right|^{\beta'}\right] \leq K u^{\be'}\Delta_n^{\frac{\beta-\beta'}{\beta}},
	\end{align*}
	using parts (c) and (g) of Lemma \ref{lemaux} as well as (\ref{ineqchi3}) which holds for any $0 < p < 2$. The claim now follows from the same reasoning as in \eqref{key:Nn_above}, with an additional step of successive conditioning.
\end{proof}

\begin{lemma}
	We have
	\begin{align*}
	\frac{1}{n-k_n-2}\mathbb{E}\left[|R_2^n(u_n)|\right]\leq K ( k_n^{-\beta/2p+\iota} \vee u\Delta_n^{1/2})
	\end{align*}
	where the constant $K$ might depend on $p$, $\be$ and $\iota$ but not on $u$.
\end{lemma}
\begin{proof}
We get 
\begin{align*}
	&\frac{1}{n-k_n-2}\E\left[\sum_{i=k_n+3}^{\Nn}\left|r_i^2(u) \1_{\Ci}\right|\right] \leq K k_n^{-\beta/2p+\iota}
	\end{align*}
with the same arguments that led to \eqref{key:Nn_above}. Similar arguments as in the previous proof plus Assumption \ref{ass:strong}, boundedness of $\Delta_n^{-p/\beta}\widetilde{V}_i^n(p)$ on $(\Ci)^C$ and $\be > 1$ to ensure the existence of moments give 
	\begin{align*}
	&\E_{i-2}^n\left[\left|r_i^2(u)\1_{(\Ci)^C}\right|\right]\\
	=&\E_{i-2}^n\left[\left|\cos\left(u\frac{\sigma_{\tau_{i-2}^n}(\widetilde{\Delta_{i}^nS}-\widetilde{\Delta_{i-1}^nS})}{\widetilde{V}_i^n(p)^{1/p}}\right)-\cos\left(u\frac{\sigma_{\tau_{i-2}^n}\left(\widetilde{\Delta_{i}^nS}-\left(\frac{\lambda_{\tau_{i-2}^n}}{\lambda_{\tau_{i-3}^n}}\right)^{\frac{1}{\beta}-1}\widetilde{\Delta_{i-1}^nS}\right)}{\widetilde{V}_i^n(p)^{1/p}}\right)\right|\1_{(\Ci)^C}\right]\\
	\leq& K u\left|\frac{\sigma_{\tau_{i-2}^n}}{\Delta_n^{-1/\beta}\widetilde{V}_i^n(p)^{1/p}}\right|\E_{i-2}^n\left[\left|\Delta_n^{-1/\beta}\widetilde{\Delta_{i-1}^nS}-\left(\frac{\lambda_{\tau_{i-2}^n}}{\lambda_{\tau_{i-3}^n}}\right)^{\frac{1}{\beta}-1}\Delta_n^{-1/\beta}\widetilde{\Delta_{i-1}^nS}\right|\1_{(\Ci)^C}\right]\\
	\leq& Ku \left| \left|\lambda_{\tau_{i-3}^n}\right|^{\frac{1}{\beta}-1}-\left|\lambda_{\tau_{i-2}^n}\right|^{\frac{1}{\beta}-1} \right|. 
	\end{align*}
	The expectation of the right hand side is bounded by $K u \De_n^{1/2}$,
	using (\ref{ineqsila}) and successive conditioning. We then obtain  
	\begin{align*}
	&\frac{1}{n-k_n-2}\E\left[\sum_{i=k_n+3}^{\Nn}\left|r_i^2(u) \1_{(\Ci)^C} \right|\right] \leq K u \Delta_n^{1/2}
	\end{align*}
	as in the previous proof. 
\end{proof}

\begin{lemma}
	We have
	\begin{align*}
	\frac{1}{n-k_n-2}\mathbb{E}\left[|R_3^n(u)|\right]\leq K (u^\beta\alpha_n \vee u^\beta k_n^{-1/2} \vee k_n^{-\beta/2p+\iota})
	\end{align*}
	where the constant $K$ might depend on $p$, $\be$ and $\iota$ but not on $u$.
\end{lemma}
\begin{proof}
	Again we obtain
\begin{align*}
	&\frac{1}{n-k_n-2}\E\left[\sum_{i=k_n+3}^{\Nn}\left|r_i^3(u) \1_{\Ci}\right|\right] \leq K k_n^{-\beta/2p+\iota}
	\end{align*}
as in \eqref{key:Nn_above}, this time using the boundedness of $x \mapsto \exp(-x)$ on the positive halfline. We then use a first order Taylor expansion of the (random) function
	\begin{align} \label{def:fi}
	f^n_{i,u}(x)=\exp\left(-\frac{A_{\beta}u^\beta|\sigma_{\tau_{i-2}^n}|^\beta|\lambda_{\tau_{j-2}^n}|^{1-\beta}((\phi_i^n)^{1-\beta}+(\phi_{i-1}^n)^{1-\beta})}{x^{\beta/p}}\right) \1_{(\Ci)^C}
	\end{align}
	(defined for $x > 0$) and get
	\begin{align*}
	&\E_{i-2}^n\left[\exp\left(-\frac{A_\beta u^\beta|\sigma_{\tau_{i-2}^n}|^\beta|\lambda_{\tau_{j-2}^n}|^{1-\beta} ((\phi_i^n)^{1-\beta}+(\phi_{i-1}^n)^{1-\beta})}{\Delta_n^{-1}\widetilde{V}_i^n(p)^{\beta/p}}\right)\right] \1_{(\Ci)^C}\\
	&~~~-\E_{i-2}^n\left[\exp\left(-\frac{C_{p,\beta} u^\beta|\sigma_{\tau_{i-2}^n}|^\beta|\lambda_{\tau_{j-2}^n}|^{1-\beta} ((\phi_i^n)^{1-\beta}+(\phi_{i-1}^n)^{1-\beta})}{(\sigs_i)^{\beta/p}}\right)\right] \1_{(\Ci)^C}\\
	=&\left(\Delta_n^{-p/\beta}\widetilde{V}_i^n(p)-\sigs_i\mu_{p,\beta}^{p/\beta}\kappa_{p,\beta}^{p/\beta}\right)\E_{i-2}\left[(f^{n}_{i,u})'(\epsilon_i^n)\right] \1_{(\Ci)^C}
	\end{align*} 
	for some $\fF_{\tau_{i-2}^n}$-measurable $\epsilon_i^n$ between $\Delta_n^{-p/\beta}\widetilde{V}_i^n(p)$ and $\sigs_i\mu_{p,\beta}^{p/\beta}\kappa_{p,\beta}^{p/\beta}$. An easy computation proves 
	\begin{align} \label{ineqfiu}
	\E_{i-2}^n\left[|(f^n_{i,u})'(X)|\right] \leq K \frac{ u^\beta}{X^{p/\beta+1}} \quad \text{and} \quad \E_{i-2}^n\left[|(f^n_{i,u})''(X)|\right]\leq K \frac{ u^\beta}{X^{p/\beta+2}}
	\end{align}
	for any positive $\fF_{\tau_{i-2}^n}$-measurable random variable $X$ where the constant $K$ does not depend on $u$. Using the boundedness of $\Delta_n^{-p/\beta}\widetilde{V}_i^n(p)$ and $\sigs_i\mu_{p,\beta}^{p/\beta}\kappa_{p,\beta}^{p/\beta}$ again we get
	\begin{align*}
	\E\left[|r_i^3(u_n)\1_{(\Ci)^C}|\right] &\leq K u_n^\beta \E\left[\left|\Delta_n^{-p/\beta}\widetilde{V}_i^n(p)-\sigs_i\mu_{p,\beta}^{p/\beta}\kappa_{p,\beta}^{p/\beta}\right|\right]\\
	&\leq K u_n^\beta \E\left[\left|\Delta_n^{-p/\beta}(\widetilde{V}_i^n(p)-\overline{V}_i^n(p))\right|+\left|\Delta_n^{-p/\beta}\overline{V}_i^n(p)-\sigs_i\mu_{p,\beta}^{p/\beta}\kappa_{p,\beta}^{p/\beta}\right|\right]\\
	&\leq K u_n^\beta (k_n^{-1/2}\vee \alpha_n \vee \Delta_n^{1/2}) \leq K u_n^\beta (k_n^{-1/2}\vee \alpha_n )
	\end{align*}
	where the last line holds by Lemma \ref{lemoverline} and \eqref{ineqoverline} plus the definition of $\al_n$.
\end{proof}

\begin{lemma}
	We have
	\begin{align*}
	\frac{1}{n-k_n-2}\mathbb{E}\left[|R_4^n(u)|\right]\leq K u^\beta\Delta_n k_n
	\end{align*}
	where the constant $K$ might depend on $p$ and $\be$ but not on $u$.
\end{lemma}
\begin{proof}
	Recall the function $f^n_{i,u}$ from (\ref{def:fi}) and set $\widetilde{r}_{i,n}=(\sigs_i-|\sigma_{\tau_{i-2}^n}|^p|\lambda_{\tau_{i-2}^n}|^{\frac{p}{\beta}-p})$. Then 
	\begin{align}
	&\E\left[\left|R_4^n(u)\right|\right] \nonumber\\ \leq&
	\E\left[\left|R_4^n(u)-\mu_{p,\beta}^{p/\beta}\kappa_{p,\beta}^{p/\beta}\sum_{i=k_n+3}^{\Nn}\E_{i-2}^n\left[(f^n_{i, u})'(\mu_{p,\beta}^{p/\beta}\kappa_{p,\beta}^{p/\beta}|\sigma_{\tau_{i-2}^n}|^p|\lambda_{\tau_{i-2}^n}|^{\frac{p}{\beta}-p})\right]\widetilde{r}_{i,n} \right|\right]\label{ineq:R41}\\
	+&\E\left[\left|\mu_{p,\beta}^{p/\beta}\kappa_{p,\beta}^{p/\beta}\sum_{i=k_n+3}^{\Nn}(\widetilde{r}_{i,n}-\E_{i-k_n-3}^n\left[\widetilde{r}_{i,n}\right])\E_{i-2}^n\left[(f^n_{i, u})'(\mu_{p,\beta}^{p/\beta}\kappa_{p,\beta}^{p/\beta}|\sigma_{\tau_{i-2}^n}|^p|\lambda_{\tau_{i-2}^n}|^{\frac{p}{\beta}-p})\right]\right|\right]\label{ineq:R42}\\
	+&\E\left[\left|\mu_{p,\beta}^{p/\beta}\kappa_{p,\beta}^{p/\beta}\sum_{i=k_n+3}^{\Nn}\E_{i-2}^n\left[(f^n_{i, u})'(\mu_{p,\beta}^{p/\beta}\kappa_{p,\beta}^{p/\beta}|\sigma_{\tau_{i-2}^n}|^p|\lambda_{\tau_{i-2}^n}|^{\frac{p}{\beta}-p})\right]\E_{i-k_n-3}^n\left[\widetilde{r}_{i,n}\right]\right|\right].\label{ineq:R43}
	\end{align}
	In the sequel we prove the same rate of convergence for all three terms on the right hand side. Starting with \eqref{ineq:R41}, from the definition of $r_i^4(u)$ we have
	\begin{align*}
	r_i^4(u)=\E_{i-2}^n\left[f_{i,u}(\mu_{p,\beta}^{p/\beta}\kappa_{p,\beta}^{p/\beta}\sigs_i)-f_{i,u}(\mu_{p,\beta}^{p/\beta}\kappa_{p,\beta}^{p/\beta}|\sigma_{\tau_{i-2}^n}|^p|\lambda_{\tau_{i-2}^n}|^{\frac{p}{\beta}-p})\right],
	\end{align*}
	using the independence of $\phi_{i-1}^n$ and $\phi_{i}^n$ from $\fF_{\tau_{i-2}^n}$. A second order Taylor expansion, possible by the usual boundedness assumptions, now gives
	\begin{align*}
	&\E\left[\left|R_4^n(u)-\sum_{i=k_n+3}^{\Nn} \mu_{p,\beta}^{p/\beta}\kappa_{p,\beta}^{p/\beta} \widetilde{r}_{i,n} \E_{i-2}^n\left[(f^n_{i, u})'(\mu_{p,\beta}^{p/\beta}\kappa_{p,\beta}^{p/\beta}|\sigma_{\tau_{i-2}^n}|^p|\lambda_{\tau_{i-2}^n}|^{\frac{p}{\beta}-p})\right]\right|\right]\\
	=& \E\left[\left|\sum_{i=k_n+3}^{\Nn}\frac12(\mu_{p,\beta}^{p/\beta}\kappa_{p,\beta}^{p/\beta})^2(\widetilde{r}_{i,n})^2\E_{i-2}^n\left[f_{i,u_n}''(\delta_{i,n})\right]\right|\right]
	\end{align*}
for some $\delta_{i,n}$ between $\mu_{p,\beta}^{p/\beta}\kappa_{p,\beta}^{p/\beta}\sigs_i$ and  $\mu_{p,\beta}^{p/\beta}\kappa_{p,\beta}^{p/\beta}|\sigma_{\tau_{i-2}^n}|^p|\lambda_{\tau_{i-2}^n}|^{\frac{p}{\beta}-p}$. Now it is an easy consequence of (\ref{ineqfiu}) and Lemma \ref{lem:ineqsila} together with the reasoning from \eqref{key:Nn_above} that the expectation in \eqref{ineq:R41} is bounded by $K u^\beta k_n$ with $K$ as in the statement of the lemma.

For \eqref{ineq:R43}, boundedness of all processes involved gives $$\left|\E_{i-2}^n\left[(f^n_{i, u})'(\mu_{p,\beta}^{p/\beta}\kappa_{p,\beta}^{p/\beta}|\sigma_{\tau_{i-2}^n}|^p|\lambda_{\tau_{i-2}^n}|^{\frac{p}{\beta}-p})\right]\right|\leq Ku^\beta$$ by (\ref{ineqfiu}) whereas Lemma \ref{lem:ineqsila}  proves
	\begin{align}\label{r_i4_bound}
		\left|\E_{i-k_n-3}\left[\widetilde{r}_{i,n}\right]\right|\leq Kk_n\Delta_n.
	\end{align}
	Thus
	\begin{align*}
	\E\left[\mu_{p,\beta}^{p/\beta}\kappa_{p,\beta}^{p/\beta}\sum_{i=k_n+3}^{\Nn}\left|\E_{i-2}^n\left[(f_{i, u})'(\mu_{p,\beta}^{p/\beta}\kappa_{p,\beta}^{p/\beta}|\sigma_{\tau_{i-2}^n}|^p|\lambda_{\tau_{i-2}^n}|^{\frac{p}{\beta}-p})\right]\right|\left|\E_{i-k_n-3}^n\left[\widetilde{r}_{i,n}\right]\right|\right]
	\leq K u^\beta k_n.
	\end{align*}

Finally, for the treatment of \eqref{ineq:R42} we have to be a little more specific. A simple computation proves 
\[
\E_{i-2}^n\left[(f^n_{i, u})'(\mu_{p,\beta}^{p/\beta}\kappa_{p,\beta}^{p/\beta}|\sigma_{\tau_{i-2}^n}|^p|\lambda_{\tau_{i-2}^n}|^{\frac{p}{\beta}-p})\right] = K u^\be L(p,u,\be) |\sigma_{\tau_{i-2}^n}|^{-p}|\lambda_{\tau_{i-2}^n}|^{p-\frac{p}{\beta}}
\]
for some $K$ as above. Thus, setting $\Xi_i=\widetilde{r}_{i,n}-\E_{i-k_n-3}\left[\widetilde{r}_{i,n}\right]$ we can bound \eqref{ineq:R42} by
\begin{align} 
& K u^\be L(p,u,\be) \E\left[\left|\sum_{i=k_n+3}^{\Nn} \Xi_i  \left(|\sigma_{\tau_{i-2}^n}|^{-p}|\lambda_{\tau_{i-2}^n}|^{p-\frac{p}{\beta}} - |\sigma_{\tau_{i-k_n-3}^n}|^{-p}|\lambda_{\tau_{i-k_n-3}^n}|^{p-\frac{p}{\beta}}\right)\right|\right] \label{ineq:R42a} \\ +& 
K u^\be L(p,u,\be) \E\left[\left|\sum_{i=k_n+3}^{\Nn} \Xi_i  |\sigma_{\tau_{i-k_n-3}^n}|^{-p}|\lambda_{\tau_{i-k_n-3}^n}|^{p-\frac{p}{\beta}}\right|\right]. \label{ineq:R42b} 
\end{align} 
An application of the Cauchy-Schwarz inequality bounds (\ref{ineq:R42a}) by the product of
\[
 K u^\be L(p,u,\be) \E\left[\sum_{i=k_n+3}^{\Nn} \Xi_i^2 \right]^{1/2}
\]
and 
\[
 \E\left[\sum_{i=k_n+3}^{\Nn} \left(|\sigma_{\tau_{i-2}^n}|^{-p}|\lambda_{\tau_{i-2}^n}|^{p-\frac{p}{\beta}} - |\sigma_{\tau_{i-k_n-3}^n}|^{-p}|\lambda_{\tau_{i-k_n-3}^n}|^{p-\frac{p}{\beta}}\right)^2\right]^{1/2}.
\]
Lemma \ref{lem:ineqsila} together with \eqref{r_i4_bound} proves 
	\begin{align}\label{eq:R4_1}
	\E_{i-k_n-3}^n[|\Xi_i|^2]\leq K k_n\Delta_n,
	\end{align}
and from Lemma \ref{lem:ineq_power} we have 
\[
\E_{i-k_n-3}^n\left[\left(|\sigma_{\tau_{i-2}^n}|^{-p}|\lambda_{\tau_{i-2}^n}|^{p-\frac{p}{\beta}} - |\sigma_{\tau_{i-k_n-3}^n}|^{-p}|\lambda_{\tau_{i-k_n-3}^n}|^{p-\frac{p}{\beta}}\right)^2\right] \le K k_n \De_n
\]
with the same reasoning as when establishing (\ref{ineqsila}). Boundedness of $L(p,u,\be)$ and $\Nn \le C \De_n^{-1}$ now prove that (\ref{ineq:R42a}) is bounded by $K u^\beta k_n$. 

For (\ref{ineq:R42b}) we use an argument involving discrete martingales, and we first change the upper summation bound from $\Nn$ to $\Nn + 2k_n +5$ as the corresponding error term is of the order $K u^\beta k_n$ by boundedness of $\si$ and $\la$, so similar to the one from  (\ref{ineq:R42a}). The martingale argument is explained the easiest if we first pretend that the factors $K u^\be L(p,u,\be) |\sigma_{\tau_{i-k_n-3}^n}|^{-p}|\lambda_{\tau_{i-k_n-3}^n}|^{p-\frac{p}{\beta}}$ were not present. We write 
\begin{align} \label{resid}
\sum_{i=k_n+3}^{\Nn+2k_n+5}\Xi_i = \sum_{j=1}^{k_n+3}A_j+ \sum_{i=\lfloor\Nn/(k_n+3)\rfloor)(k_n+3) +2k_n +6}^{\Nn+2k_n +5}\Xi_i
\end{align}
with 
\begin{align*}
A_j=\sum_{i=1}^{\lfloor\Nn /(k_n+3)\rfloor +1}\Xi_{i(k_n+3)+(j-1)} =\sum_{i=1}^{\infty}\Xi_{i(k_n+3)+(j-1)} \1_{\{i-1 \le \lfloor\Nn /(k_n+3)\rfloor \}},
\end{align*}
$j=1, \ldots, k_n+3.$ It can be shown that
\begin{align*}
&\E[\Xi_{i(k_n+3)+(j-1)} \1_{\{i-1 \le \lfloor\Nn /(k_n+3)\rfloor \}} 
\Xi_{\ell(k_n+3)+(j-1)} \1_{\{\ell-1 \le \lfloor\Nn /(k_n+3)\rfloor \}}] = 0
\end{align*}
holds for every $1 \le \ell < i$, using the fact that by construction one knows at time $\tau_{(i-1)(k_n+3)+(j-1)}^n$ whether the event $\{\tau_{(i-1)(k_n+3)}^n \le 1 \}$ has happened or not. The latter event is equivalent to $\{i-1 \le \lfloor\Nn/(k_n+3)\rfloor \}$, so after conditioning on $\fF_{(i-1)(k_n+3)+(j-1)}^n$ the claim follows from $\E_{r-k_n-3}^n[\Xi_{r}] = 0$ for every $r$.

Thus, by \eqref{eq:R4_1}, Cauchy-Schwarz inequality and $\lfloor\Nn /(k_n+3)\rfloor +1 \leq Kn/k_n$ we obtain
	\begin{align*}
	\E\left[\left|A_j\right|\right]\leq K\E\left[\left|\sum_{i=1}^{\lfloor\Nn /(k_n+3)\rfloor +1}|\Xi_{k_n+3+(j-1)+(i-1)(k_n+1)}|^2\right|\right]^{1/2} \le K.
	\end{align*}
As the sum over the residual terms in (\ref{resid}) has at most $k_n+2$ elements we obtain 
	\begin{align*}
	\E\left[\left|\sum_{i=k_n+3}^{\Nn+2k_n+5}\Xi_i \right|\right]\leq K k_n 
	\end{align*}
as desired. If we now include $K u^\be L(p,u,\be) |\sigma_{\tau_{i-k_n-3}^n}|^{-p}|\lambda_{\tau_{i-k_n-3}^n}|^{p-\frac{p}{\beta}}$, we just get an additional factor $u^{\be}$ as usual. This is due to boundedness of $\si$ and $\la$ again (and of the function $L$) plus the fact that measurability w.r.t.\ $\fF_{\tau_{i-k_n-3}^n}$ keeps the martingale property from above intact.
\end{proof}

\begin{lemma}\label{lemma:Z}
	We have
	\begin{align*}
	&\frac{1}{n-k_n-2}\E\left[|Z^n(u)-\overline{Z}^n(u)|\right]\\ \leq& K \left(k_n^{-\beta/2p+\iota}\vee \Delta_n^{1/2}(u^{\beta/2-\iota} + u^{\beta/2})\left(k_n^{-1/2} \vee \alpha_n \vee (k_n\Delta_n)^{1/2}\right)^{1/2}\right)
	\end{align*}
where the constant $K$ might depend on $p$, $\be$ and $\iota$ but not on $u$.
\end{lemma}

\begin{proof}
As usual we have 
\begin{align*}
	&\frac{1}{n-k_n-2}\E\left[\sum_{i=k_n+3}^{\Nn}\left|(z_i^n(u) - \overline z_i^n(u)) \1_{\Ci}\right|\right] \leq K k_n^{-\beta/2p+\iota},
	\end{align*}
and for the analogous sum involving $\1_{(\Ci)^C}$, as in the previous proof, we may change the upper summation index to $\Nn+2$ without loss of generality. Now, note 
that by Lemma \ref{lem:exp} and using the same arguments for $z_i(u_n)$ we have
	\begin{align*}
	\E_{i-2}^n\left[z_i(u)\right]=0=\E_{i-2}^n\left[\overline{z}_i(u)\right].
	\end{align*}
Thus for all $i,j \ge k_n +3$ with $j-i \geq 2$ we have 
	\begin{align*}
	\E\left[(z_i(u)-\overline{z}_i(u))(z_j(u)-\overline{z}_j(u))\1_{(\Ci)^C}\1_{(\mathcal{C}_j^n)^C}\1_{\{\Nn+2\geq j\}}\1_{\{\Nn+2\geq i\}}\right]=0
	\end{align*}
where we have used that $\1_{\{\Nn+2\geq j\}},\1_{\{\Nn+2\geq i\}},\Ind_{(\Ci)^C}$ and $\1_{(\mathcal{C}_j^n)^C}$ are all $\mathcal{F}_{\tau_{j-2}^n}^n$-measurable.
Using $2|xy|\leq x^2+y^2$ and $\Nn \le C \De_n^{-1}$ we then get
	\begin{align*}
	&\E \left[\left(\sum_{i=k_n+3}^{\Nn+2}(z_i(u)-\overline{z}_i(u))\1_{(\Ci)^C}\right)^2\right]\\
	=&\E\bigg[\sum_{i=k_n+3}^{\infty}(z_i(u)-\overline{z}_i(u))\1_{(\mathcal{C}_i^n)^C}\1_{\{\Nn+2\geq i\}}\big((z_i(u)-\overline{z}_i(u))\1_{(\mathcal{C}_i^n)^C}\1_{\{\Nn+2\geq i\}}\\
	+&(z_{i-1}(u)-\overline{z}_{i-1}(u))\1_{(\mathcal{C}_{i-1}^n)^C}\1_{\{\Nn+2\geq i-1\}}+(z_{i+1}(u)-\overline{z}_{i+1}(u))\1_{(\mathcal{C}_{i+1}^n)^C}\1_{\{\Nn+2\geq i+1\}}\big)\bigg]\\
	 \leq& 3 \sum_{i=k_n+3}^{\lfloor C \De_n^{-1} \rfloor +2} \E\left[(z_i(u)-\overline{z}_i(u))^2\1_{(\Ci)^C}\right].
	\end{align*}

Now, with (\ref{identL}) plus the standard inequalities $|\cos(x)-\cos(y)|^2\leq4|x-y|^p$ and $|\exp(-x)-\exp(-y)|^2\leq|x-y|^p$, which hold for all $p \in (0,2]$, we obtain 
	\begin{align*}
	&\E_{i-2}^n\left[|(z_i(u)-\overline{z}_i(u))\1_{(\Ci)^C}|^2\right] \nonumber\\ 
	\leq&  2\E_{i-2}^n\left[\left|\cos\left(u\frac{\sigma_{\tau_{i-2}^n}\left(\widetilde{\Delta_{i}^nS}-\left(\frac{\lambda_{\tau_{i-2}^n}}{\lambda_{\tau_{i-3}^n}}\right)^{\frac{1}{\beta}-1}\widetilde{\Delta_{i-1}^nS}\right)}{\widetilde{V}_i^n(p)^{1/p}}\right)\right.\right.\nonumber\\
	&~~~~~~~~~~~~~~~\left.\left.-\cos\left(u\frac{\lambda_{\tau_{i-2}^n}^{1-1/\beta}\widetilde{\Delta_i^n S}-\lambda_{\tau_{i-3}^n}^{1-1/\beta}\widetilde{\Delta_{i-1}^nS}}{\Delta_n^{1/\beta}\mu_{p,\beta}^{1/\beta}\kappa_{p,\beta}^{1/\beta}}\right)\right|^2\1_{(\Ci)^C}\right]  \nonumber\\ 
	+& 2\E_{i-2}^n\left[\left|\E_{i-2}^n\left[\exp\left(-\frac{A_\beta u^\beta{|\sigma_{\tau_{i-2}^n}|}^\beta|\lambda_{\tau_{i-2}^n}|^{1-\beta} ((\phi_i^n)^{1-\beta}+(\phi_{i-1}^n)^{1-\beta})}{\Delta_n^{-1}\widetilde{V}_i^n(p)^{\beta/p}}\right)\right]\right.\right.\nonumber\\ &~~~~~~~~~~~~~~~-\left.\left.\E_{i-2}^n\left[\exp(-u^\beta C_{p,\beta}((\phi_{i}^n)^{1-\beta}+(\phi_{i-1}^n)^{1-\beta})\vphantom{\frac{1}{2}}\right]\right|^2\1_{(\Ci)^C}\right] \nonumber\\
	 \leq& Ku^{\beta-\iota}\E_{i-2}^n\left[\left|\frac{|\sigma_{\tau_{i-2}^n}||\lambda_{\tau_{i-2}^n}|^{\frac{1}{\beta}-1}}{\widetilde{V}_i^n(p)^{1/p}}-\frac{1}{\Delta_n^{1/\beta}\mu_{p,\beta}^{1/\beta}\kappa_{p,\beta}^{1/\beta}}\right|^{\beta-\iota}\1_{(\Ci)^C}|\lambda_{\tau_{i-2}^n}^{1-1/\beta}\widetilde{\Delta_i^n S}-\lambda_{\tau_{i-3}^n}^{1-1/\beta}\widetilde{\Delta_{i-1}^nS}|^{\beta-\iota}\right]\nonumber\\
	& \qquad + Ku^\beta\E_{i-2}^n\left[\left|\frac{A_\beta |\sigma_{\tau_{i-2}^n}|^\beta|\lambda_{\tau_{i-2}^n}|^{1-\beta}}{\Delta_n^{-1}\widetilde{V}_i^n(p)^{\beta/p}}-\frac{A_\beta}{\mu_{p,\beta}\kappa_{p,\beta}}\right|\1_{(\Ci)^C}|(\phi_{i}^n)^{1-\beta}+(\phi_{i-1}^n)^{1-\beta}| \right], 
	\end{align*}
	and part (a) of Lemma \ref{lemaux} together with $\E\left|(\phi_{i}^n)^{1-\beta}+(\phi_{i-1}^n)^{1-\beta}\right|<\infty$ and the $\fF_{i-2}^n$-measurability of the other terms proves that the term above is bounded by
	\begin{align*}
	&K u^{\beta-\iota} \left|\frac{|\sigma_{\tau_{i-2}^n}||\lambda_{\tau_{i-2}^n}|^{\frac{1}{\beta}-1}\mu_{p,\beta}^{1/\beta}\kappa_{p,\beta}^{1/\beta}-\Delta_n^{-1/\beta}\widetilde{V}_i^n(p)^{1/p}}{\Delta_n^{-1/\beta}\widetilde{V}_i^n(p)^{1/p}\mu_{p,\beta}^{1/\beta}\kappa_{p,\beta}^{1/\beta}}\right|^{\beta-\iota} \1_{(\Ci)^C}
	\\
	&\quad + Ku^{\beta} \left|\frac{|\sigma_{\tau_{i-2}^n}|^\beta|\lambda_{\tau_{i-2}^n}|^{1-\beta}\mu_{p,\beta}\kappa_{p,\beta}-\Delta_n^{-1}\widetilde{V}_i^n(p)^{\beta/p}}{\Delta_n^{-1}\widetilde{V}_i^n(p)^{\beta/p}\mu_{p,\beta}\kappa_{p,\beta}}\right|\1_{(\Ci)^C}.
	\end{align*}
	Let us for a moment only discuss the first term. On $(\Ci)^C$ and using Condition \ref{ass:strong}, $\Delta_n^{-p/\beta}\widetilde{V}_i^n(p)$ as well as all quantities involving $\si$ and $\la$ are bounded from above and below by $K$. Thus, together with $|x^q-y^q| \leq q|\max(x,y)^{q-1}||x-y|$ for $q\geq1$ we get
	\begin{align*}
	&\left|\frac{|\sigma_{\tau_{i-2}^n}||\lambda_{\tau_{i-2}^n}|^{\frac{1}{\beta}-1}\mu_{p,\beta}^{1/\beta}\kappa_{p,\beta}^{1/\beta}-\Delta_n^{-1/\beta}\widetilde{V}_i^n(p)^{1/p}}{\Delta_n^{-1/\beta}\widetilde{V}_i^n(p)^{1/p}\mu_{p,\beta}^{1/\beta}\kappa_{p,\beta}^{1/\beta}}\right|^{\beta-\iota} \Ind_{(\Ci)^C}\\
	&\le K \left|\frac{1}{p}\max\left(|\sigma_{\tau_{i-2}^n}|^p|\lambda_{\tau_{i-2}^n}|^{\frac{p}{\beta}-p}\mu_{p,\beta}^{p/\beta}\kappa_{p,\beta}^{p/\beta},\Delta_n^{-p/\beta}\widetilde{V}_i^n(p)\right)^{1/p-1}\Ind_{(\Ci)^C}\right.\\
	&~~~~~~~~~~~~~~\left.\left| \Delta_n^{-p/\beta}\widetilde{V}_i^n(p)-|\sigma_{\tau_{i-2}^n}|^p|\lambda_{\tau_{i-2}^n}|^{\frac{p}{\beta}-p}\mu_{p,\beta}^{p/\beta}\kappa_{p,\beta}^{p/\beta}\right|\right|^{\beta-\iota}\\
	&\leq K \left|\Delta_n^{-p/\beta}\widetilde{V}_i^n(p)-|\sigma_{\tau_{i-2}^n}|^p|\lambda_{\tau_{i-2}^n}|^{\frac{p}{\beta}-p}\mu_{p,\beta}^{p/\beta}\kappa_{p,\beta}^{p/\beta}\right|^{\beta-\iota} \Ind_{(\Ci)^C}.
	\end{align*}
With a similar argument for the second term we then obtain 
	\begin{align*}
		&\E_{i-2}^n\left[|(z_i(u)-\overline{z}_i(u))\1_{(\Ci)^C}|^2\right]  \\ &\leq  \1_{(\Ci)^C} K u^{\beta-\iota} \left|\Delta_n^{-p/\beta}\widetilde{V}_i^n(p)-|\sigma_{\tau_{i-2}^n}|^p|\lambda_{\tau_{i-2}^n}|^{\frac{p}{\beta}-p}\mu_{p,\beta}^{p/\beta}\kappa_{p,\beta}^{p/\beta}\right|^{\beta-\iota} \\ &+
 \1_{(\Ci)^C} K u^{\beta} \left|\Delta_n^{-p/\beta}\widetilde{V}_i^n(p)-|\sigma_{\tau_{i-2}^n}|^p|\lambda_{\tau_{i-2}^n}|^{\frac{p}{\beta}-p}\mu_{p,\beta}^{p/\beta}\kappa_{p,\beta}^{p/\beta}\right|.
	\end{align*}
Now, we have 
\begin{align*}
&\E\left[\left|\Delta_n^{-p/\beta}\widetilde{V}_i^n(p)-|\sigma_{\tau_{i-2}^n}|^p|\lambda_{\tau_{i-2}^n}|^{\frac{p}{\beta}-p}\mu_{p,\beta}^{p/\beta}\kappa_{p,\beta}^{p/\beta}\right|\right] \\ &\leq  \E\left[\left|\Delta_n^{-p/\beta}\widetilde{V}_i^n(p)-\Delta_n^{-p/\beta}\overline{V}_i^n(p)\right|+\left|\Delta_n^{-p/\beta}\overline{V}_i^n(p)-\sigs_i\mu_{p,\beta}^{p/\beta}\kappa_{p,\beta}^{p/\beta}\right|\right.\\
	&\quad\quad+\left.\left||\sigma_{\tau_{i-2}^n}|^p|\lambda_{\tau_{i-2}^n}|^{\frac{p}{\beta}-p}\mu_{p,\beta}^{p/\beta}\kappa_{p,\beta}^{p/\beta}-\sigs_i\mu_{p,\beta}^{p/\beta}\kappa_{p,\beta}^{p/\beta}\right|\vphantom{\frac{1}{2}}\right]\\
	&\leq \left(k_n^{-1/2} \vee \alpha_n \vee \Delta_n^{1/2} \vee (k_n\Delta_n)^{1/2}\right)
\end{align*}
using Lemma \ref{lemoverline}, (\ref{ineqas}) and Lemma \ref{lem:ineqsila}. A similar result holds with the exponent being replaced by $\beta - \iota > 1$. The claim now follows easily. 
\end{proof}

\subsection{Proof of the main theorems}

\subsubsection{Proof of Theorem \ref{thmlncon}}
As discussed before, we may assume Conditions \ref{ass:strong} and \ref{ass:strongstopping} to hold. First, we have
\[
\widetilde{L}^n(p,u) -  L(p,u,\beta) = \frac{1}{\Nn-k_n-2}\sum_{i=k_n+3}^{\Nn} \left(\cos\left(u\frac{\widetilde{\Delta_i^n X}-\widetilde{\Delta_{i-1}^nX} }{(\widetilde{V}_i^n(p))^{1/p}}\right) - L(p,u,\beta) \right),
\]
and it is a simple consequence of (\ref{conv:Nn}), the decomposition in (\ref{eq:decomp}) and Lemmas \ref{lemma:R1} to \ref{lemma:Z} that
\[
\widetilde{L}^n(p,u) -  L(p,u,\beta) = \frac{1}{\Nn-k_n-2} \overline Z^{n}(u) + o_\P(1)
\]
holds. Note that we have convergence to zero of all bounds in Lemmas \ref{lemma:R1} to \ref{lemma:Z} as $\De_n \to 0$, $k_n \to \infty$ and $k_n \De_n \to 0$ by assumption.

The proof of   
\[
\frac{1}{\Nn-k_n-2} \overline Z^{n}(u) \pn 0
\]
follows along the lines of Lemma \ref{lemma:Z}. We may first change the upper summation index to $\Nn + 2$ which does not change anything asymptotically because of (\ref{conv:Nn}), and we then have
\[
\E\left[\overline{z}_i(u) \overline{z}_j(u) \1_{\{\Nn+2\geq j\}}\1_{\{\Nn+2\geq i\}}\right]=0
\]
for all $i,j \ge k_n +3$ with $j-i \geq 2$ using Lemma \ref{lem:exp}. 
Thus, we obtain
\[
\E\left[\left(\sum_{i=k_n+3}^{\Nn+2}\overline{z}_i(u) \right)^2\right] \leq 3 \sum_{i=k_n+3}^{\lfloor C \De_n^{-1} \rfloor +2} \E\left[(\overline{z}_i(u))^2\right] \le K \De_n^{-1},
\]
and the claim follows from (\ref{conv:Nn}) again. \qed

\subsubsection{Proof of Remark \ref{be2}}
From the definition of $\hat \beta(p, u_n, v_n)$ it follows easily that the claim $\hat \beta(p, u_n, v_n) \le 2$ is equivalent to 
	\begin{align}
 \sum_{i=k_n+3}^{\Nn}\rho^2(1-\cos\left(u_na_i\right)) \leq \sum_{i=k_n+3}^{\Nn}(1-\cos\left(\rho u_na_i\right))\label{sufficient_bound}
	\end{align}
	where we have used the shorthand notation $a_i=\frac{\widetilde{\Delta_i^n X}-\widetilde{\Delta_{i-1}^nX} }{(\widetilde{V}_i^n(p))^{1/p}}$. A sufficient condition for 
\eqref{sufficient_bound} is $\rho^2(1-\cos(x)) \leq 1-\cos(\rho x)$ for all $x\in\R$ 	which itself is equivalent to
	\begin{align}\label{sufficient_bound_2}
	g_\rho(x)= 1 - \cos\left(\rho x\right) - \rho^2(1-\cos(x)) \geq 0 \quad \text{for all }x\in\R. 
	\end{align}
	Using properties of the cosine and inserting $\rho=1/2$ we note $g_{\frac12}(x)=g_{\frac12}(-x)$ and $g_{\frac12}(x) = g_{\frac12}(x+4\pi)$. For \eqref{sufficient_bound_2} to hold it then suffices to show $g_{\frac12}(x)\geq 0$ for all $x\in[0,2\pi]$. So let $x\in[0,2\pi]$. Then
	\begin{align*}
	g_{\frac12}'(x) &= \frac{1}{2}\sin\left(\frac x2\right)- \frac14\sin(x)\\
	&=\frac12\sin\left(\frac x2\right) - \frac12 \sin\left(\frac x2\right)\cos\left(\frac x2\right) = \frac12\sin\left(\frac x2\right)\left(1-\cos\left(\frac x2\right)\right) \geq 0
	\end{align*}
	by properties of the trigonometric functions. The claim follows from $g_{\frac12}(0) = 0$.\qed

\subsubsection{Proof of Theorem \ref{thmlnclt}}
We will assume throughout that Assumptions \ref{ass:strong} and \ref{ass:strongstopping} are in place, and we set $k_n\sim C_1 \De_n^{-\varpi}$ for some $C_1> 0$ and some $\varpi\in(0,1)$ as well as $u_n\sim C_2 \De_n^{\varrho}$ for some $C_2 > 0$ and $\varrho \in (0,1)$. 
\begin{lemma} \label{lem:cond}
Under the conditions 
\begin{align*}
	\beta'<\frac{\beta}{2}, \quad \frac{1}{3}\vee\frac{1}{8\varrho}<p<\frac{\beta}{2}, \quad \varpi\geq\frac{2}{3}, \quad \frac{1}{3\beta}<\varrho< \frac{1}{\beta}, \quad \frac{1}{\beta}<\frac{\varpi}{p}-\varrho, \quad 2\varpi-\varrho\beta<1, 
	\end{align*}
we have 
\begin{align*}	
\frac{\sqrt{N_n(1)}}{u_n^{\beta/2}}\left((\widetilde{L}^n(p,u_n)-L(p,u_n,\beta)) - \frac{1}{\Nn-k_n-2} \overline Z^n(u_n)\right)  \pn 0.
\end{align*}
\end{lemma}

\begin{proof}
The claim follows from a tedious but straightforward computation, using (\ref{conv:Nn}) and Lemmas \ref{lemma:R1} to \ref{lemma:Z} as well as the conditions on $k_n$ and $u_n$. 
\end{proof}

\begin{lemma}\label{Lem:calcvar}
	Let $u_n$ be as above and set $v_n=\rho u_n$ for some $\rho > 0$. Then, for every fixed $i > 2$ we have 
	\begin{align*}
	\frac{1}{u_n^{\beta/2}v_n^{\beta/2}}\E\left[\overline{z}_i(u_n)\overline{z}_i(v_n)\right]&\to C_{p,\beta}\kappa_{\beta,\beta}\frac{2+2\rho^\beta-|1-\rho|^\beta-(1+\rho)^\beta}{2\rho^{\beta/2}},\\
	\frac{1}{u_n^{\beta/2}v_n^{\beta/2}}\E\left[\overline{z}_i(u_n)\overline{z}_{i-1}(v_n)\right]&\to C_{p,\beta}\kappa_{\beta,\beta}\frac{2+2\rho^\beta-|1-\rho|^\beta-(1+\rho)^\beta}{4\rho^{\beta/2}},
	\end{align*}
	and the same result holds with interchanged roles of $u_n$ and $v_n$. 
\end{lemma}
\begin{proof}
Using the shorthand notation $\widehat{\Delta_i^nS} = \lambda_{\tau_{i-2}^n}^{-1/\beta+1}\widetilde{\Delta_{i}^nS}$ and the equality $\cos(x)\cos(y)=\frac{1}{2}\left(\cos(x-y)+\cos(x+y)\right)$ we have 
		\begin{align*}
		&\cos\left(u_n\frac{\widehat{\Delta_i^nS}-\widehat{\Delta_{i-1}^nS}}{\Delta_n^{1/\beta}\mu_{p,\beta}^{1/\beta}\kappa_{p,\beta}^{1/\beta}}\right)\cos\left(v_n\frac{\widehat{\Delta_i^nS}-\widehat{\Delta_{i-1}^nS}}{\Delta_n^{1/\beta}\mu_{p,\beta}^{1/\beta}\kappa_{p,\beta}^{1/\beta}}\right)\\
		=&\frac12 \left(\cos\left(\frac{(u_n-v_n)\widehat{\Delta_i^nS}+(-u_n+v_n)\widehat{\Delta_{i-1}^nS}}{\Delta_n^{1/\beta}\mu_{p,\beta}^{1/\beta}\kappa_{p,\beta}^{1/\beta}}\right)+\cos\left(\frac{(u_n+v_n)\widehat{\Delta_i^nS}+(-u_n-v_n)\widehat{\Delta_{i-1}^nS}}{\Delta_n^{1/\beta}\mu_{p,\beta}^{1/\beta}\kappa_{p,\beta}^{1/\beta}}\right)\right),\\
		&\cos\left(u_n\frac{\widehat{\Delta_i^nS}-\widehat{\Delta_{i-1}^nS}}{\Delta_n^{1/\beta}\mu_{p,\beta}^{1/\beta}\kappa_{p,\beta}^{1/\beta}}\right)\cos\left(v_n\frac{\widehat{\Delta_{i-1}^nS}-\widehat{\Delta_{i-2}^nS}}{\Delta_n^{1/\beta}\mu_{p,\beta}^{1/\beta}\kappa_{p,\beta}^{1/\beta}}\right)\\
		=&\frac12 \left(\cos\left(\frac{u_n\widehat{\Delta_i^nS}+(-u_n-v_n)\widehat{\Delta_{i-1}^nS}+v_n\widehat{\Delta_{i-2}^nS}}{\Delta_n^{1/\beta}\mu_{p,\beta}^{1/\beta}\kappa_{p,\beta}^{1/\beta}}\right)\right.\\
		&~~~~~~~~~~~~~~~~~~~~~~~~~~~~\left.+\cos\left(\frac{u_n\widehat{\Delta_i^nS}+(-u_n+v_n)\widehat{\Delta_{i-1}^nS}-v_n\widehat{\Delta_{i-2}^nS}}{\Delta_n^{1/\beta}\mu_{p,\beta}^{1/\beta}\kappa_{p,\beta}^{1/\beta}}\right)\right).
		\end{align*}
	With the same notation as in the proof of Lemma \ref{lem:exp} we obtain
	\begin{align*}
	&(u_n-v_n)\Delta_n^{-1/\beta}\widehat{\Delta_{i}^nS}+(-u_n+v_n)\Delta_n^{-1/\beta}\widehat{\Delta_{i-1}^nS}\nonumber\\
	\sim& u_n(1-\rho)\lambda_{\tau_{i-2}^n}^{-1/\beta+1}((\phi_i^n\lambda_{\tau_{i-2}^n})^{1-\beta})^{1/\beta}S^{(1)}+u_n(\rho-1)\lambda_{\tau_{i-3}^n}^{-1/\beta+1}((\phi_{i-1}^n\lambda_{\tau_{i-3}^n})^{1-\beta})^{1/\beta}S^{(2)}\nonumber\\
	\sim& u_n|1-\rho|S^{(3)}((\phi_i^n)^{1-\beta}+(\phi_{i-1}^n)^{1-\beta})^{1/\beta}
	\end{align*}
	as the ($\mathcal{F}_{\tau_{i-2}^n}$-conditional) distribution,	and in the same manner
	\begin{align*}
	(u_n+v_n)\Delta_n^{-1/\beta}\widehat{\Delta_{i}^nS}+&(-u_n-v_n)\Delta_n^{-1/\beta}\widehat{\Delta_{i-1}^nS}  \sim u_n(1+\rho)S^{(3)}((\phi_i^n)^{1-\beta}+(\phi_{i-1}^n)^{1-\beta})^{1/\beta}\nonumber\\
	u_n\Delta_n^{-1/\beta}\widehat{\Delta_{i}^nS}+&(-u_n-v_n)\Delta_n^{-1/\beta}\widehat{\Delta_{i-1}^nS}+v_n\Delta_n^{-1/\beta}\widehat{\Delta_{i-2}^nS}\nonumber\\ &\sim u_nS^{(3)}((\phi_i^n)^{1-\beta}+(1+\rho)^\beta (\phi_{i-1}^n)^{1-\beta}+ \rho^\beta (\phi_{i-2}^n)^{1-\beta})^{1/\beta},\nonumber\\
	u_n\Delta_n^{-1/\beta}\widehat{\Delta_{i}^nS}+&(-u_n+v_n)\Delta_n^{-1/\beta}\widehat{\Delta_{i-1}^nS}-v_n\Delta_n^{-1/\beta}\widehat{\Delta_{i-2}^nS}\nonumber\\ &\sim u_nS^{(3)}((\phi_i^n)^{1-\beta}+|1-\rho|^\beta (\phi_{i-1}^n)^{1-\beta}+ \rho^\beta (\phi_{i-2}^n)^{1-\beta})^{1/\beta}.
	\end{align*}
	We see in particular that exchanging the roles of $u_n$ and $v_n$ is irrelevant to the distributions. Then, with Lemma \ref{lem:exp} and its proof,
	\begin{align*}
	&\frac{1}{u_n^{\beta/2}v_n^{\beta/2}}\E\left[\overline{z}_i(u_n)\overline{z}_i(v_n)\right]\\
	=&\frac{1}{2u_n^{\beta}\rho^{\beta/2}}\E\left[\exp\left(-C_{p,\beta}u_n^\beta|1-\rho|^\beta ((\phi_i^n)^{1-\beta}+(\phi_{i-1}^n)^{1-\beta})\right)\right.\\
	&\quad\left.+\exp\left(-C_{p,\beta}u_n^\beta(1+\rho)^\beta ((\phi_i^n)^{1-\beta}+(\phi_{i-1}^n)^{1-\beta})\right)\right]\\
	&\quad-\frac{1}{u_n^{\beta}\rho^{\beta/2}}\E\left[\exp\left(-u_n^\beta C_{p,\beta}((\phi_i^n)^{1-\beta}+(\phi_{i-1}^n)^{1-\beta})\right)\right]\E\left[\exp\left(-u_n^\beta\rho^\beta C_{p,\beta}((\phi_i^n)^{1-\beta}+(\phi_{i-1}^n)^{1-\beta})\right)\right],\\
	&\frac{1}{u_n^{\beta/2}v_n^{\beta/2}}\E\left[\overline{z}_i(u_n)\overline{z}_{i-1}(v_n)\right]\\
	=&\frac{1}{2u_n^{\beta/2}v_n^{\beta/2}}\E\left[\exp\left(-C_{p,\beta}u_n^\beta((\phi_i^n)^{1-\beta}+(1+\rho)^\beta (\phi_{i-1}^n)^{1-\beta}+ \rho^\beta (\phi_{i-2}^n)^{1-\beta})\right)\right.\\
	&\quad\left.+\exp\left(-C_{p,\beta}u_n^\beta((\phi_i^n)^{1-\beta}+|1-\rho|^\beta (\phi_{i-1}^n)^{1-\beta}+ \rho^\beta (\phi_{i-2}^n)^{1-\beta})\right)\right]\\
	&\quad-\frac{1}{u_n^{\beta/2}v_n^{\beta/2}}\E\left[\exp\left(-u_n^\beta C_{p,\beta}((\phi_i^n)^{1-\beta}+(\phi_{i-1}^n)^{1-\beta})\right)\right]\E\left[\exp\left(-u_n^\beta\rho^\beta C_{p,\beta}((\phi_i^n)^{1-\beta}+(\phi_{i-1}^n)^{1-\beta})\right)\right].
	\end{align*}
	We now have for example
	\begin{align*}
	& \E\left[\exp\left(-C_{p,\beta}u_n^\beta |1-\rho|^\beta ((\phi_i^n)^{1-\beta}+(\phi_{i-1}^n)^{1-\beta})\right)\right]\\ =&-C_{p,\beta}u_n^\beta |1-\rho|^\beta\E\left[\exp(-\epsilon_{i}^n)((\phi_i^n)^{1-\beta}+(\phi_{i-1}^n)^{1-\beta})\right]+1
	\end{align*}
	for some $\epsilon_{i}^n\in[0,u_n^\beta |1-\rho|^\beta C_{p,\beta}((\phi_i^n)^{1-\beta}+(\phi_{i-1}^n)^{1-\beta})]$, and since $u_n \to 0$ and $\E[(\phi_i^n)^{1-\beta}] = M < \infty$ hold, dominated convergence gives
	\begin{align*}
	\E\left[\exp(-\epsilon_i^n)((\phi_i^n)^{1-\beta}+(\phi_{i-1}^n)^{1-\beta})\right] \rightarrow \E\left[((\phi_i^n)^{1-\beta}+(\phi_{i-1}^n)^{1-\beta})\right]=\kappa_{\beta,\beta}.
	\end{align*}
The same argument for the other terms, including the usage of $u_n \to 0$ when dealing with the product of the two expectations, gives
	\begin{align*}
	&\frac{1}{u_n^{\beta/2}v_n^{\beta/2}}\E\left[\overline{z}_i(u_n)\overline{z}_i(v_n)\right]\\
	\to& \frac{-C_{p,\beta}|1-\rho|^\beta \kappa_{\beta,\beta}+1-C_{p,\beta}(1+\rho)^\beta \kappa_{\beta,\beta}+1-2\left(- C_{p,\beta}\kappa_{\beta,\beta}- \rho^\beta C_{p,\beta}\kappa_{\beta,\beta}+1\right)}{2\rho^{\beta/2}}\nonumber\\
	=&C_{p,\beta}\kappa_{\beta,\beta}\frac{2+2\rho^\beta-|1-\rho|^\beta-(1+\rho)^\beta}{2\rho^{\beta/2}}.
	\end{align*}
	Similar arguments lead to
	\begin{align*}
	\frac{1}{u_n^{\beta/2}v_n^{\beta/2}}\E\left[\overline{z}_i(u_n)\overline{z}_{i-1}(v_n)\right]&\pn C_{p,\beta}\frac{\kappa_{\beta,\beta}}{2}\frac{4+4\rho^\beta-(1+(1+\rho)^\beta+\rho^\beta)-(1+|1-\rho|^\beta+\rho^\beta)}{2\rho^{\beta/2}}\\
	&=C_{p,\beta}\kappa_{\beta,\beta}\frac{2+2\rho^\beta-|1-\rho|^\beta-(1+\rho)^\beta}{4\rho^{\beta/2}}.
	\end{align*}
\end{proof}

\begin{lemma}\label{lem:Zn_conv}
Let $u_n$ be as above and set $v_n=\rho u_n$ for some $\rho > 0$. Then the $\fF$-stable convergence in law
\begin{align*}
\left(\frac{\sqrt{\De_n}}{u_n^{\beta/2}}\overline{Z}^n(u_n),\frac{\sqrt{\De_n}}{v_n^{\beta/2}}\overline{Z}^n(v_n)\right)\tols (X,Y)
\end{align*}
holds, where $(X,Y)$ is mixed normal distributed with mean $0$ and covariance matrix $\mathcal{C}$ consisting of
	\begin{align*}
	&\mathcal{C}_{11}=\mathcal{C}_{22}=\int_{0}^{1} \frac{1}{\lambda_s}ds~C_{p,\beta}\kappa_{\beta,\beta}(4-2^\beta),\\
	&\mathcal{C}_{12}=\mathcal{C}_{21}=\int_{0}^{1} \frac{1}{\lambda_s}ds~C_{p,\beta}\kappa_{\beta,\beta}\frac{2+2\rho^\beta - (1+\rho)^\beta  -|1-\rho|^\beta }{\rho^{\beta/2}}.
	\end{align*}
\end{lemma}

\begin{proof}
We set
	\begin{footnotesize}
		\begin{align*}
		\zeta_i^n&=\frac{\sqrt{\De_n}}{u_n^{\beta/2}}\left(\overline{z}_i(u_n),\overline{z}_i(v_n)\right)\\
		&=\left(\frac{\sqrt{\De_n}}{u_n^{\beta/2}}\left(\cos\left(u_n\frac{\widehat{\Delta_{i}^nS}-\widehat{\Delta_{i-1}^nS}}{\Delta_n^{1/\beta}\mu_{p,\beta}^{1/\beta}\kappa_{p,\beta}^{1/\beta}}\right)-L(p,u_n,\beta)\right),\frac{\sqrt{\De_n}}{v_n^{\beta/2}}\left(\cos\left(v_n\frac{\widehat{\Delta_{i}^nS}-\widehat{\Delta_{i-1}^nS}}{\Delta_n^{1/\beta}\mu_{p,\beta}^{1/\beta}\kappa_{p,\beta}^{1/\beta}}\right)-L(p,v_n,\beta)\right)\right),
		\end{align*}
	\end{footnotesize}
	
\noindent and we can see by the uniform boundedness of $\overline{z}_i(\cdot)$ 
	and because of $\De_nu_n^{-\beta} \to 0$ that
	\begin{align*}
	\sum_{i=k_n+3}^{\Nn}\zeta_i^n \quad \text{and} \quad \sum_{i=k_n+3}^{\Nn+1}\left(\zeta_i^n-\E^n_{i-1}\left[\zeta_i^n\right]+\E^n_{i}\left[\zeta_{i+1}^n\right]\right)
	\end{align*}
	are asymptotically equivalent. We note that $\Nn+1$ is a $(\mathcal{F}_{\tau_{i}^n})_{i\geq1}$-stopping time and therefore in order to apply Theorem 2.2.15 in \cite{discret} it is sufficient to show that for $q=\frac{2}{1-\varrho\beta}+2>2$, $\eta_i^n=\zeta_i^n-\E^n_{i-1}\left[\zeta_i^n\right]+\E^n_{i}\left[\zeta_{i+1}^n\right]$
	\begin{align}
	&\sum_{i=k_n+3}^{\Nn+1}\E^n_{i-1}\left[\eta_i^n\right]\xrightarrow{\mathbb{P}}(0,0),\label{conv:1}\\
	&\sum_{i=k_n+3}^{\Nn+1}\left(\E^n_{i-1}\left[\eta_{i}^{n,j}\eta_{i}^{n,k}\right]-\E^n_{i-1}\left[\eta_{i}^{n,j}\right]\E^n_{i-1}\left[\eta_{i}^{n,k}\right]\right)\xrightarrow{\mathbb{P}}\mathcal{C}_{jk},\label{conv:2}\\
	&\sum_{i=k_n+3}^{\Nn+1}\E^n_{i-1}\left[\lVert\zeta_i^n\rVert^q\right] \xrightarrow{\mathbb{P}}0, \label{conv:3}\\
	&\label{eq:stableconv}
\sum_{i=k_n+3}^{\Nn+1}\E^n_{i-1}\left[\zeta_i^n(M_{\tau_{i}^n}-M_{\tau_{i-1}^n})\right]\xrightarrow{\mathbb{P}}0,
	\end{align} 
hold, where $M$ is either one of the Brownian motions $W$, $\widetilde{W}$ or $\overline W$ or a bounded martingale orthogonal to any of the Brownian motions. 

First note that Lemma \ref{lem:exp} gives $\E^n_{i-1}\left[\zeta_{i+1}^n\right]=(0,0)$ and therefore $\E^n_{i-1}\left[\eta_i^n\right]=
(0,0)$ by definition as well. \eqref{conv:1} then holds. Also, $\varrho<\frac{1}{\beta}$ gives
	\begin{align*}
	\frac{1}{1-\varrho\beta} - \varrho\beta\frac{1+1-\varrho\beta}{1-\varrho\beta}=\frac{(1-\varrho\beta)^2}{1-\varrho\beta}>0.
	\end{align*} 
Thus, the uniform boundedness of $\overline{z}_i(\cdot)$ and Assumption \ref{ass:strongstopping} give 
	\begin{align*}
	\sum_{i=k_n+3}^{\Nn+1}\E^n_{i-1}\left[\left|\sqrt{\frac{\De_n}{u_n^{\beta}}}\overline{z}_i(u_n)\right|^q\right]\leq \Nn u_n^{-\beta\frac{1+1-\varrho\beta}{1-\varrho\beta}}\De_n^{1+\frac{1}{1-\varrho\beta}}\leq C u_n^{-\beta\frac{1+1-\varrho\beta}{1-\varrho\beta}}\De_n^{\frac{1}{1-\varrho\beta}} \to 0
	\end{align*}
	which proves \eqref{conv:3}. To show \eqref{conv:2} we first recall $\E^n_{i-1}\left[\eta_i^n\right]=(0,0)$ and then a simple calculation yields
	\begin{align*}
	&\E^n_{i-1}\left[\left(\zeta_i^{n,j}-\E^n_{i-1}\left[\zeta_i^{n,j}\right]+\E^n_{i}\left[\zeta_{i+1}^{n,j}\right]\right)\left(\zeta_i^{n,k}-\E^n_{i-1}\left[\zeta_i^{n,k}\right]+\E^n_{i}\left[\zeta_{i+1}^{n,k}\right]\right)\right]\\
	&=\E^n_{i-1}\left[\zeta_i^{n,j}\zeta_i^{n,k}\right]-\E^n_{i-1}\left[\zeta_i^{n,j}\right]\E^n_{i-1}\left[\zeta_i^{n,k}\right]+\E^n_{i-1}\left[\zeta_i^{n,j}\zeta_{i+1}^{n,k}\right]+\E^n_{i-1}\left[\zeta_i^{n,k}\zeta_{i+1}^{n,j}\right]\\
	&\qquad+\E\left[\E^n_{i}\left[\zeta_{i+1}^{n,j}\right]\E^n_{i}\left[\zeta_{i+1}^{n,k}\right]\right],
	\end{align*}
	using iterated expectations, $\E^n_{i-1}\left[\zeta_{i+1}^n\right]=(0,0)$ and the fact that the distribution of $\E^n_{i}\left[\zeta_{i+1}^{n,j}\right]\E^n_{i}\left[\zeta_{i+1}^{n,k}\right]$ is independent of $\mathcal{F}_{\tau_{i-1}^n}$. We then prove 
	\begin{align*}
	\sum_{i=k_n+3}^{\Nn}\E^n_{i-1}\left[\zeta_i^{n,j}\zeta_i^{n,k}\right]&\pn \int_{0}^{1} \frac{1}{\lambda_s}ds\lim_{n\rightarrow\infty} \De_n^{-1}\E\left[\zeta_m^{n,j}\zeta_m^{n,k}\right],\\
	\sum_{i=k_n+3}^{\Nn}\E^n_{i-1}\left[\zeta_i^{n,j}\right]\E^n_{i-1}\left[\zeta_i^{n,k}\right]&\pn \int_{0}^{1} \frac{1}{\lambda_s}ds\lim_{n\rightarrow\infty} \De_n^{-1} \E\left[\E^n_{m-1}\left[\zeta_m^{n,j}\right]\E^n_{m-1}\left[\zeta_m^{n,k}\right]\right],\\
	\sum_{i=k_n+3}^{\Nn}\E^n_{i-1}\left[\zeta_i^{n,j}\zeta_{i+1}^{n,k}\right]&\pn \int_{0}^{1} \frac{1}{\lambda_s}ds\lim_{n\rightarrow\infty} \De_n^{-1} \E\left[\zeta_m^{n,j}\zeta_{m+1}^{n,k}\right],\\
	\sum_{i=k_n+3}^{\Nn}\E^n_{i-1}\left[\zeta_i^{n,k}\zeta_{i+1}^{n,j}\right]&\pn \int_{0}^{1} \frac{1}{\lambda_s}ds\lim_{n\rightarrow\infty} \De_n^{-1} \E\left[\zeta_m^{n,k}\zeta_{m+1}^{n,j}\right]
	\end{align*}
	and that the limits on the right-hand side exist and are the same, irrespective of the choice of $m$. This would result in
	\begin{align}\label{zgwb:varconv}
	\sum_{i=k_n+3}^{\Nn}&\left(\zeta_i^n-\E^n_{i-1}\left[\zeta_i^n\right]+\E^n_{i}\left[\zeta_{i+1}^n\right]\right)^T\left(\zeta_i^n-\E^n_{i-1}\left[\zeta_i^n\right]+\E^n_{i}\left[\zeta_{i+1}^n\right]\right)\\
	&\xrightarrow{\mathbb{P}}\int_{0}^{1} \frac{1}{\lambda_s}ds\lim_{n\rightarrow\infty}n\left(\E\left[({\zeta_m^n})^T\zeta_m^n\right]+\E\left[(\zeta_m^n)^T\zeta_{m+1}^n\right]+\E\left[(\zeta_{m+1}^n)^T\zeta_{m}^n\right]\right),\nonumber
	\end{align}
	everything for an arbitrary $m$.
	
	We give the arguments for the first convergence result above in detail, the other ones can be treated in exactly the same way. We set $X_i^{n,1}=\E_{i}^n\left[({\zeta_{i+1}^n})^T({\zeta_{i+1}^n})\right]$ and then prove 
	\begin{align} \label{eq:convev}
	\sum_{i=k_n+2}^{\Nn+1}\E_{i-1}^n\left[X_i^{n,1}\right]\pn \int_{0}^{1} \frac{1}{\lambda_s}ds\lim_{n\rightarrow\infty} \De_n^{-1} \E\left[X_m^{n,1}\right]
	\end{align}
	as well as
	\begin{align} \label{eq:asympneg}
	\sum_{i=k_n+2}^{\Nn+1}\E_{i-1}^n\left[((X_i^{n,1})_{jk})^2\right]\pn 0, \quad j,k=1,2.
	\end{align}
	Lemma 2.2.12 in \cite{discret} (plus the usual asymptotic negligibility when adding finitely many summands) finally gives the claim. 
	Now, note first that the distribution of $X_i^{n,1}$ is independent of $\mathcal{F}_{\tau_{i-1}^n}$. Therefore
	\begin{align*}
	\E_{i-1}^n\left[X_i^{n,1}\right]=\E\left[X_i^{n,1}\right] \quad \text{and} \quad \E_{i-1}^n\left[(X_i^{n,1})^2\right]=\E\left[(X_i^{n,1})^2\right]. 
	\end{align*} 
	(\ref{eq:convev}) is then an easy consequence of 
	\begin{align*}
\sum_{i=k_n+2}^{\Nn+1}\E_{i-1}^n\left[X_i^{n,1}\right]&=\De_n(\Nn-k_n) \De_n^{-1}\E\left[X_i^{n,1}\right] \pn \int_{0}^{1} \frac{1}{\lambda_s}ds\lim_{n\rightarrow\infty} \De_n^{-1}\E\left[X_m^{n,1}\right]
\end{align*}
where we used \eqref{conv:Nn} and Lemma \ref{Lem:calcvar} to prove that the limit on the right hand side exists. To show (\ref{eq:asympneg}) we use
Jensen inequality to obtain
\begin{align*}
\E\left[((X_i^{n,1})_{jk})^2\right]
&\leq \E\left[\E_i^n\left[\left(\zeta_{i+1}^{n,j}\right)^2\right]\right]\left\|\left(\zeta_{i+1}^{n,k}\right)^2\right\|_\infty,
\end{align*}
and uniform boundedness of $\overline{z}_i(\cdot)$ as well as $\Nn \le C \De_n^{-1}$ give
\begin{align*}
\sum_{i=k_n+2}^{\Nn+1}\E\left[((X_i^{n,1})_{jk})^2\right]&\leq K\sum_{i=k_n+2}^{\Nn+1}\frac{\De_n}{u_n^{\beta}}\E\left[\left(\zeta_{i+1}^{n,j}\right)^2\right]
\to 0,
\end{align*}
using $\De_nu_n^{-\beta}\to 0$ along with Lemma \ref{Lem:calcvar} in the last step. The proof of (\ref{conv:2}) can then be finished by a tedious but straightforward computation, combining \eqref{zgwb:varconv} with Lemma \ref{Lem:calcvar} again.

Finally, to prove \eqref{eq:stableconv} we use Theorem 4.34 in Chapter III of \cite{jacshi2003}. We set for $k_n+3\leq i \leq \Nn$ and $t\geq \tau_{i-2}^n$:
\begin{align*}
\mathcal{H}:=\mathcal{F}_{\tau_{i-2}^n} \quad \text{and} \quad \mathcal{H}_t := \mathcal{H} \bigvee \sigma\left(S_r : t \ge r\geq \tau_{i-2}^n \right),
\end{align*} 
i.e. $(\mathcal{H}_t)_{t\geq \tau_{i-2}^n }$ is the filtration generated by $\mathcal{H}$ and $\sigma\left(S_r : t \ge r \geq \tau_{i-2}^n \right)$. Now $(S_t)_{t\geq \tau_{i-2}^n }$ is a process with independent increments w.r.t.\ to $\sigma\left(S_r : r\geq \tau_{i-2}^n \right)$. For all $t\geq \tau_{i-2}^n$ we set $K_t:= \E\left[\zeta_i|\mathcal{H}_t\right]$ and note that $K_{\tau_{i}^n} = \zeta_i$ due to $\zeta_i$ being $\mathcal{H}_{\tau_i^n}$-measurable. Then with the aforementioned Theorem 4.34 we have
\begin{align*}
\zeta_i = K_{\tau_{i}^n} =  K_{\tau_{i-2}^n} +\int_{\tau_{i-2}}^{\tau_{i}}H_s dS_s,
\end{align*}
where $(H_t)_{t\geq \tau_{i-2}^n}$ is a predictable process. Then
\begin{align*}
\E_{i-1}^n\left[\zeta_i^n(M_{\tau_{i}}-M_{\tau_{i-1}})\right] &=
\left(K_{\tau_{i-2}^n}+\int_{\tau_{i-2}}^{\tau_{i-1}}H_sdS_s\right) \E_{i-1}^n\left[M_{\tau_{i}}-M_{\tau_{i-1}}\right]\\  &~~~~~~~~~~~~~~~~~~~~+\E_{i-1}^n\left[\int_{\tau_{i-1}}^{\tau_{i}}H_sdS_s(M_{\tau_{i}}-M_{\tau_{i-1}})\right]= 0
\end{align*}
where we used that the martingale $(S_t)_{t\geq 0}$ is orthogonal to $M$ in all cases.
\end{proof}

\begin{corollary} \label{cor520}
Suppose that the conditions in Lemma \ref{lem:cond} hold and choose $v_n=\rho u_n$ for some $\rho > 0$. Then we have the $\fF$-stable convergence in law 
\begin{align*}
\left(\frac{\sqrt{N_n(1)}}{u_n^{\beta/2}}(\widetilde{L}^n(p,u_n)-L(p,u_n,\beta)),\frac{\sqrt{N_n(1)}}{v_n^{\beta/2}}(\widetilde{L}^n(p,v_n)-L(p,v_n,\beta))\right)\tols (X',Y')
\end{align*}
where $(X',Y')$ is jointly normal distributed with mean $0$ and covariance matrix $\mathcal{C'}$ consisting of
	\begin{align*}
	\mathcal{C}_{11}'=\mathcal{C}_{22}'=C_{p,\beta}\kappa_{\beta,\beta}(4-2^\beta), \quad \mathcal{C}_{12}'=\mathcal{C}_{21}'=C_{p,\beta}\kappa_{\beta,\beta}\frac{2+2\rho^\beta - (1+\rho)^\beta  -|1-\rho|^\beta }{\rho^{\beta/2}}.
	\end{align*}
\end{corollary}

\begin{proof}
Using Lemma \ref{lem:cond} we have
\[
\frac{\sqrt{N_n(1)}}{u_n^{\beta/2}}\left((\widetilde{L}^n(p,u_n)-L(p,u_n,\beta)) - \frac{1}{\Nn-k_n-2} \overline Z^n(u_n)\right)  \pn 0,
\]
and similarly for $v_n$, and from several applications of (\ref{conv:Nn}) together with $k_n \De_n \to 0$ we also get
\[
\frac{\sqrt{N_n(1)}}{\sqrt{\De_n}(\Nn-k_n-2)} = \frac{1}{\sqrt{\De_n \Nn}} (1+o_\P(1)) \pn \frac 1{\sqrt{\int_{0}^{1}\frac{1}{\lambda_s}ds}}. 
\]
Then Lemma \ref{lem:Zn_conv} together with the properties of stable convergence in law yields the claim. 
\end{proof}

\subsubsection{Proof of Theorem \ref{thm:beta_hat}}
Using $L(p,u_n,\beta) = \E \left[\exp \left(-u^\beta C_{p,\beta}((\phi^{(1)})^{1-\beta}+(\phi^{(2)})^{1-\beta})\right) \right]$, a Taylor expansion of the function 
\begin{align*}
	(x,y)\mapsto\frac{\log(-(x-1))-\log(-(y-1))}{\log(1/\rho)} 
	\end{align*} 
	with gradient 
	\begin{align*}
	(g_1(x),g_2(y))=\left(\frac{1}{\log(1/\rho)(x-1)},\frac{1}{\log(1/\rho)(1-y)}\right)
	\end{align*} 
	around $(L(p,u_n,\beta),L(p,v_n,\beta))$ gives
			\begin{align}
	u_n^{\beta/2}&\sqrt{N_n(1)}(\hat{\beta}(p,u_n,v_n)-\beta)= \nonumber\\
	&u_n^{\beta/2}\sqrt{N_n(1)}\left(\frac{\log(-(L(p,u_n,\beta)-1))-\log(-(L(p,v_n,\beta)-1))}{\log(u_n/v_n)}-\beta\right) \label{thm2:bias}\\
	+&\frac{1}{\log(u_n/v_n)}\frac{u_n^\beta}{\E[\exp(-u_n^\beta C_{p,\beta}((\phi^{(1)})^{1-\beta}+(\phi^{(2)})^{1-\beta})]-1}\frac{\sqrt{N_n(1)}}{u_n^{\beta/2}}(\widetilde{L}^n(p,u_n)-L(p,u_n,\beta))\label{thm1:normapprox}\\
	+&\frac{1}{\log(u_n/v_n)}\frac{1}{\rho^{\beta/2}}\frac{v_n^\beta}{1-\E[\exp(-v_n^\beta C_{p,\beta}((\phi^{(1)})^{1-\beta}+(\phi^{(2)})^{1-\beta})]}\frac{\sqrt{N_n(1)}}{v_n^{\beta/2}}(\widetilde{L}^n(p,v_n)-L(p,v_n,\beta))\label{thm2:normapprox}\\
	+&u_n^{\beta}(g_1(\eta_1^n)-g_1(L(p,u_n,\beta)))\frac{\sqrt{N_n(1)}}{u_n^{\beta/2}}(\widetilde{L}^n(p,u_n)-L(p,u_n,\beta)) \label{ineqeta}\\
	+&v_n^{\beta}\frac{1}{\rho^{\beta/2}}(g_2(\eta_2^n)-g_2(L(p,v_n,\beta)))\frac{\sqrt{N_n(1)}}{v_n^{\beta/2}}(\widetilde{L}^n(p,v_n)-L(p,v_n,\beta)) \nonumber,
	\end{align}
	for some $\eta_1^n$ between $\widetilde{L}^n(p,u_n)$ and $L(p,u_n,\beta)$ and some $\eta_2^n$ between $\widetilde{L}^n(p,v_n)$ and $L(p,v_n,\beta)$.

	As before, we have
	\begin{align*}
	1-L(p,u_n,\beta)=\E\left[\exp(-\epsilon_{1}^n) C_{p,\beta}u_n^\beta \left((\phi^{(1)})^{1-\beta}+(\phi^{(2)})^{1-\beta}\right)\right] 
	\end{align*}
	for some $\eps_1^n$ between $0$ and $C_{p,\beta}u_n^\beta ((\phi^{(1)})^{1-\beta}+(\phi^{(2)})^{1-\beta})$. Obviously, $\eps_1^n \to 0$ almost surely, so by dominated convergence 
	\begin{align} \label{expans}
	\frac{1-L(p,u_n,\beta)}{u_n^\be} \to \E\left[C_{p,\beta} \left((\phi^{(1)})^{1-\beta}+(\phi^{(2)})^{1-\beta}\right)\right] = C_{p,\beta} \ka_{\be, \be} > 0.  
	\end{align}
We now prove 
\begin{align} \label{expans2}
u_n^{\beta} \left |\frac 1{\eta_1^n-1} - \frac 1{L(p,u_n,\beta))-1} \right| = u_n^{\beta} \frac {|L(p,u_n,\beta)-\eta_1^n|}{(1-\eta_1^n)(1-L(p,u_n,\beta))}  \pn 0 
\end{align}
from which, together with Corollary \ref{cor520} and Slutsky's lemma, the asymptotic negligibility of (\ref{ineqeta}) follows. A similar result obviously holds for the term involving $\eta_2^n$ and $v_n$. Using (\ref{expans}), we get (\ref{expans2}) from
\[
\frac {1-\eta_1^n}{|L(p,u_n,\beta)-\eta_1^n|} \pn \infty.
\]
To prove the latter claim we use $1-\eta_1^n = (1-L(p,u_n,\be)) + (L(p,u_n,\beta)-\eta_1^n)$ and the fact that $(1-L(p,u_n,\be))$ is of the order $u_n^\be$ using (\ref{expans}) while $|L(p,u_n,\beta)-\eta_1^n| \le |L(p,u_n,\beta)-\widetilde{L}^n(p,u_n)|$ is at most of order $\De_n^{1/2} u_n^{\be/2}$ using Corollary \ref{cor520}. $\De_nu_n^{-\beta} \to 0$ together with Slutsky's lemma then yields the claim.

	In order to prove the convergence of the bias term \eqref{thm2:bias} towards zero we use that for $\eps_1^n$ as above there exists some $\epsilon_{2}^n$ between $\E[\exp(-\epsilon_{1}^n)((\phi^{(1)})^{1-\beta}+(\phi^{(2)})^{1-\beta})]$ and $\kappa_{\beta,\beta}$ such that
	\begin{align*}
	&\log(-(L(u_n,p,\beta)-1))=\log(u_n^\beta C_{p,\beta})+\log\left(\E\left[\exp(-\epsilon_{1}^n)\left((\phi^{(1)})^{1-\beta}+(\phi^{(2)})^{1-\beta}\right)\right]\right)\\
	=&\log(u_n^\beta C_{p,\beta})+\frac{1}{\epsilon_{2}^n}\left(\E[\exp(-\epsilon_{1}^n)((\phi^{(1)})^{1-\beta}+(\phi^{(2)})^{1-\beta})]-\kappa_{\beta,\beta})+\log(\kappa_{\beta,\beta}\right).
	\end{align*}
	Clearly $\eps_2^n \to \ka_{\be, \be}$, and with $\eps_3^n$ between $0$ and $\eps_1^n$ we obtain for an arbitrary  $\iota>0$ 
	\begin{align*}
	&\frac{1}{u_n^{\beta-\iota}\epsilon_{2}^n}\left(\E\left[\exp(-\epsilon_1^n)\left((\phi^{(1)})^{1-\beta}+(\phi^{(2)})^{1-\beta}\right)\right]-\kappa_{\beta,\beta}\right)\\
	=&\frac{1}{\epsilon_{2}^n}\E\left[\exp(-\epsilon_3^n)\frac{(-\epsilon_{1}^n)}{u_n^{\beta-\iota}}\left((\phi^{(1)})^{1-\beta}+(\phi^{(2)})^{1-\beta}\right)\right] \to 0,
	\end{align*}
	where we used $\epsilon_1^n ((\phi^{(1)})^{1-\beta}+(\phi^{(2)})^{1-\beta}) \leq u_n^\beta C_{p,\beta} ((\phi^{(1)})^{1-\beta}+(\phi^{(2)})^{1-\beta})^2$ and dominated convergence via part (c) of Assumption \ref{ass:stopping}. The same arguments hold for $\log(-(L(v_n,p,\beta)-1))$, so
	\begin{align*}
	&\frac{1}{u_n^{\beta-\iota}}\left(\frac{\log(-(L(p,u_n,\beta)-1))-\log(-(L(p,v_n,\beta)-1))}{\log(u_n/v_n)}-\beta\right) \\
	=&\frac{1}{u_n^{\beta-\iota}}\left(\frac{\log(u_n^\beta C_{p,\beta}) +\log(\kappa_{\beta,\beta})-(\log(v_n^\beta C_{p,\beta})+\log(\kappa_{\beta,\beta}))}{\log(u_n/v_n)}-\beta\right) + o_p(1) \pn 0.
	\end{align*}
	Now, $\frac{1}{3\beta}<\varrho$ and $N_n(1)\leq C \De_n^{-1}$ yield $u_n^{\frac{3}{2}\beta-\iota}\sqrt{N_n(1)}\rightarrow 0$ almost surely for $\iota > 0$ small enough, and therefore \eqref{thm2:bias} converges in probability to zero. 
	
	The claim then follows from deriving the asymptotics of (\ref{thm1:normapprox}) and \eqref{thm2:normapprox} for which we use (\ref{expans}) and Corollary \ref{cor520}. The form of the limiting variance is then computed easily.
\qed

\subsubsection{Proof of Corollary \ref{cor:beta_hat}}
The result follows from Theorem \ref{thm:beta_hat} immediately, upon using Slutsky's lemma and $u_n^{\beta/2}\sqrt{N_n(1)} \to \infty$ almost surely, where the latter is a consequence of (\ref{conv:Nn}) and $\De_n u_n^{-\be} \to 0$. \qed

\bibliographystyle{chicago}
\bibliography{literatur}
\end{document}